\documentclass[12pt]{amsart}
\usepackage[english]{babel}
\usepackage{tikz}
\usetikzlibrary{calc,arrows,patterns}
\usepackage[utf8]{inputenc}
     \newcommand{\PARENS}[1]{\left(#1\right)}
          \newcommand{\ccases}[1]{\begin{cases}#1\end{cases}}
         \usepackage{amsmath}
        \usepackage{mathrsfs}
        \usepackage{amssymb,amsthm,amscd}
       \usepackage{upgreek}
\newcommand{\coloneq}{:=}

\usepackage{bigstrut}
\usepackage{array}
\usepackage{mathtools,arydshln}\mathtoolsset{centercolon}
\usepackage[
pdftex,hyperfootnotes]{hyperref}
\hypersetup{pdfauthor={Gerardo Araznibarreta and Manuel Ma\~{n}as}, pdftitle={Multivariate Laurent orthogonal polynomials and integrable systems}, pdfsubject={Mathematics, Mathematical Physics},pdfkeywords={Multivariate orthogonal polynomials, Borel--Gauss factorization, quasi-determinants, Christoffel-Darboux kernels, Darboux transformations,  Christoffel formula, integrable hierarchies, Toda equations, KP equations}, pdfcreator={PDF.creator@UCM},pdfproducer={PDF.producer@UCM}}

\usepackage[total={18cm,24.5cm},top=2cm, left=1.5cm]{geometry}
\usepackage[toc,page]{appendix}
\usepackage{tocvsec2}

\newcommand{\ee}{\boldsymbol{e}}
\newcommand{\n}{\boldsymbol{n}}

\newcommand{\x}{\boldsymbol{x}}

\newcommand{\z}{\boldsymbol{z}}
\newcommand{\q}{{\boldsymbol{\alpha}}}

\renewcommand{\d}{\operatorname{d}}

\newcommand{\Z}{\mathbb{Z}}
\newcommand{\R}{\mathbb{R}}
\newcommand{\C}{\mathbb{C}}
\newcommand{\I}{\mathbb{I}}

\newtheorem{pro}{Proposition}[section]
\newtheorem{lemma}{Lemma}[section]
\newtheorem{definition}{Definition}[section]
\newtheorem{theorem}{Theorem}[section]
\newtheorem{cor}{Corollary}[section]
\numberwithin{equation}{section}

\begin{document}

	\title[Darboux transformations for multivariate orthogonal  polynomials]{Darboux transformations for \\multivariate orthogonal  polynomials}
	\author{Gerardo Ariznabarreta}\address{Departamento de F\'{\i}sica Te\'{o}rica II (M\'{e}todos Matem\'{a}ticos de la F\'{\i}sica), Universidad Complutense de Madrid, 28040-Madrid, Spain}
	\email{gariznab@ucm.es}
	\thanks{GA thanks economical support from the Universidad Complutense de Madrid  Program ``Ayudas para Becas y Contratos Complutenses Predoctorales en Espa\~{n}a 2011"}
	\author{Manuel Ma\~{n}as}
	\email{manuel.manas@ucm.es}
	\thanks{MM thanks economical support from the Spanish ``Ministerio de Econom\'{\i}a y Competitividad" research project MTM2012-36732-C03-01,  \emph{Ortogonalidad y aproximaci\'{o}n; teor\'{\i}a y aplicaciones}}
	\begin{abstract}
		Darboux transformations for  polynomial perturbations  of a real multivariate measure are found. 
		The 1D Christoffel formula is extended to the multidimensional realm:  multivariate orthogonal polynomials are expressed in terms of  last  quasi-determinants 
		and  sample matrices.  The coefficients of these matrices  are the original orthogonal polynomials evaluated at a set of nodes, which is supposed to be  poised.	A discussion for the existence of  poised sets  is given in terms of algebraic hypersufaces in the complex affine space. 
	\end{abstract}
	\maketitle
	\tableofcontents
\section{Introduction}

In a recent paper \cite{GM1} we studied how the Gauss--Borel or $LU$ factorization of  a moment matrix allows for a better understanding of the links between multivariate orthogonal polynomials (MVOPR) on a multidimensional real space $\R^D$, $D\geq 1$, and integrable systems of Toda and KP type. In particular, it was shown  how the $LU$ decomposition allows for a simple construction of the three term relation or the Christoffel--Darboux formula. Remarkably, it is also useful for the construction of Miwa type expressions in terms of quasi-tau matrices of the MVOPR or the finding  of the Darboux transformation. Indeed, we presented for the first time Darboux transformations for orthogonal polynomials in several variables, that we called elementary, and its iteration, resulting in a Christoffel quasi-determinantal type formula. These Darboux transformations allow for the construction of new MVOPR, associated with a perturbed measure, from the MVOPR of a given non perturbed measure.  Observe that they also provide a direct method to construct new solutions of the associated Toda--KP type integrable systems.

What we called elementary Darboux transformations in \cite{GM1} where given as the multiplication of the non perturbed measure  by a degree one multivariate polynomial. The $m$-th iteration of these so called  elementary Darboux transformations leads, therefore, to a perturbation by a  multivariate polynomial of degree $m$. This way of proceeding  was  motivated by the one dimensional  situation, in that case happens that the irreducible polynomials have degree one (the fundamental theorem of algebra). But, in higher dimension the situation is much richer and we do have many irreducible polynomials of higher degree.
Therefore, the territory explored for the Darboux transformations in \cite{GM1} was only a part, significant but incomplete, of a further more vast territory. In this paper, see Theorem \ref{theorem:the deal}, we give a generalization of the Darboux transformations found in \cite{GM1} that holds for  a perturbation by a polynomial of any degree.  This provides us with an elegant quasi-determinantal expression for the new MVOPR which is a broad extension of the 1D determinantal Christoffel formula.

For the construction of the mentioned general Darboux transformation we use multivariate interpolation theory, see \cite{olver}. Therefore, we need of poised sets for which the sample matrix is not singular. 
In this paper we initiate the study of poised sets  for general Darboux transformations. We find that the analysis can be splitted into two parts, one measure-independent part depending exclusively on the relative positions of nodes in the algebraic hypersurface of the generating polynomial, and another related to the non perturbed measure and the corresponding Jacobi matrices.
The geometrical part, as usual in interpolation theory,  requires of the concourse of Vandermonde matrices. In fact, of multivariate Vandermonde matrices, see \cite{olver}, or multivariate confluent Vandermonde matrices.

With the aid of basic facts in algebraic geometry, see for example \cite{hartshorne} or \cite{shafarevich} we are able to show, for generating polynomials that can be expressed as the product $\mathcal Q=\mathcal Q_1\cdots\mathcal Q_N$ of  $N$ prime factors, --see Theorem \ref{H}-- that there exists, in the complex domain,  poised sets of nodes by forbidding its belonging to any further algebraic hypersurface, different from the algebraic hypersurface of $\mathcal Q$, of certain degrees. Moreover, we see that for a perturbation of the measure by a polynomial of the form $\mathcal Q=\mathcal R^d$,  poised sets never exists, and the Darboux transformation as presented in Theorem \ref{theorem:the deal} is not applicable. However, with the use of Wro\'{n}ski matrices we can avoid this problem and find an appropriate extension of the Darboux transformation, see Theorem \ref{theorem:the big deal}, for a generating polynomial  of the form $\mathcal Q=\mathcal Q_1^{d_1}\cdots\mathcal Q_N^{d_N}$ where the polynomials $\mathcal Q_i$ are irreducible. The discussion on poised sets in this general scenario is given in Theorem \ref{HH}, where again the set of nodes when  poised can not belong to any further algebraic hypersurfaces of certain type.

The layout of the paper is as follows. Within this introduction we further perform a number observations regarding the historical background and context of the different mathematical issues discussed in this paper. Then, we reproduce, for the reader commodity, some  necessary material from \cite{GM1}. In \S \ref{S:2} we give the  Darboux transformation generated by a multivariate polynomial, and in \S \ref{S:3} we discuss poised sets, giving several conditions for the nodes in order to constitute a  poised set. Finally, in \S \ref{S:4}, we see that the previous construction fails in some cases, and them we present an extension of the Darboux transformations which overcomes this problem. 
\subsection{Historical background and context}
\subsubsection{Darboux transformations} These transformations were introduced in  \cite{darboux}  in the context of the Sturm--Liouville theory and since them have been applied in several problems. It was in \cite{matveev}, a paper where explicit solutions of the Toda lattice where found, where this covariant transformation was given the name of \emph{Darboux}. It has been used in the 1D realm of orthogonal polynomials quite successfully, see for example \cite{grunbaum-haine,yoon, bueno-marcellan1,bueno-marcellan2, marcellan}.
In Geometry,  the theory of transformations of surfaces preserving some given properties conforms a classical subject, in the list of such transformations given in the classical treatise by Einsehart \cite{eisenhart} we find the Lévy (Lucien) transformation, which later on was named as elementary Darboux transformation and known in the orthogonal polynomials context as Christoffel transformation \cite{yoon, szego}; in this paper we have denoted it by $T$. The adjoint elementary Darboux or adjoint Lévy transformation $T^{-1}$ is also relevant  \cite{matveev,dsm} and is referred some times as a Geronimus transformation \cite{yoon}. For further information see \cite{rogers-schief,gu}. 
For the iteration of elementary Darboux transformations let us mention that  Szeg\H{o} \cite{szego}  points out that for $\d\mu=\d x$  the iteration formula is due to Christoffel  \cite{christoffel}. This fact was rediscovered much  latter in the Toda context, see for example the formula (5.1.11) in \cite{matveev} for $W^+_n(N)$.

\subsubsection{Multivariate orthogonal polynomials} We refer the reader to monographs \cite{Dunkl} and \cite{xu4}.  The  recurrence relation for orthogonal polynomials in several variables was studied by Xu in \cite{xu0}, while in \cite{xu1} he linked multivariate orthogonal polynomials with a commutative family of self-adjoint operators and the spectral theorem was used to show the existence of a three term relation for the orthogonal polynomials. He discusses in \cite{xu2}  how the three term relation leads to the construction of multivariate orthogonal polynomials and cubature formul{\ae}.  Xu considers in  \cite{xu8}  polynomial subspaces that contain discrete multivariate orthogonal polynomials  with respect to the bilinear form and shows that the discrete orthogonal polynomials still satisfy a three-term relation and that Favard's theorem holds.  The analysis of orthogonal polynomials and cubature formul{\ae} on the unit ball, the standard simplex, and the unit sphere  \cite{xu6 } lead to conclude  the strong connection of orthogonal structures and cubature formul{\ae} for these three regions. The paper \cite{xu3} presents a systematic study of the common zeros of polynomials in several variables which are related to higher dimensional quadrature.
Karlin and McGregor  \cite{karlin}  and Milch \cite{milch} discussed interesting examples of multivariate Hahn and Krawt\-chouk polynomials related to growth birth and death processes. There have been substantial developments since 1975, for instance, the spectral properties of these multivariate Hahn and Krawtchouk polynomials have been studied in \cite{geronimo-iliev}. 
A  study  of two-variable orthogonal polynomials associated with a moment functional satisfying the two-variable analogue of the Pearson differential equation and an extension of some of the usual characterizations of the classical orthogonal polynomials in one variable was found \cite{fernandez}.

\subsubsection{ Quasi-determinants} For its construction we may use   Schur complements. Besides its name observe that the Schur complement was  not introduced by Issai Schur but by Emilie Haynsworth in 1968 in \cite{schur1,schur2}. In fact, Haynsworth coined that name because the Schur determinant formula given in what today is known as Schur lemma in \cite{schur}.  In the book \cite{zhang} one can find an ample overview on the Schur complement and many of its applications.  The easiest examples of quasi-determinants are Schur complements.
In the late 1920  Archibald Richardson \cite{quasideterminant1,quasideterminant2}, one of the two responsible of Littlewood--Richardson rule,  and the famous logician Arend Heyting \cite{quasideterminant3}, founder of intuitionist logic, studied possible extensions of the determinant notion to division rings. Heyting defined the \emph{designant} of a matrix with noncommutative entries, which for $2\times 2$ matrices was the Schur complement, and generalized to larger dimensions by induction. Let us stress that both Richardson's and Heyting's \emph{ quasi-determinants }
were generically rational functions of the matrix coefficients.
A definitive impulse to the modern theory was given by the Gel'fand's school \cite{quasidetermiant6,quasidetermiant7,quasidetermiant8,gelfand}.
Quasi-determinants where defined over free division rings and it was early noticed that it was not an analog of the commutative determinant but rather of a ratio determinants. A cornerstone for  quasi-determinants is  the  \emph{heredity principle}, quasi-determinants of quasi-determinants are quasi-determinants; there is no analog of such a principle for determinants.
However, many of the properties of determinants extend to this case, see the cited papers. Let us mention that in the early 1990 the Gelf'and school \cite{quasidetermiant7} already noticed the role quasi-determinants had for some integrable systems. All this paved the route, using the connection with orthogonal polynomials \emph{\`{a} la Cholesky}, to the appearance of quasi-determinants in the multivariate orthogonality context. Later, in 2006 Peter Olver applied quasi-determinants to multivariate interpolation \cite{olver}, now the blocks have different sizes, and so multiplication of blocks is only allowed if they are \emph{compatible}. In general, the (non-commutative) multiplication makes sense if the number of columns and rows of the blocks involved fit well. Moreover, we are only permitted to invert diagonal entries that in general makes the minors expansions by columns or rows not applicable but allows for other result, like the Sylvester's theorem,  to hold in this wider scenario. The last quasi-determinant used in this paper is the one described in \cite{olver}, see also \cite{GM1}.

\subsubsection{$LU$ factorization} This technique was the corner stone  for Mark Adler and Pierre van Moerbeke when  in a series of papers  where the theory
of the 2D Toda hierarchy and what they called the discrete KP hierarchy  was analyzed \cite{adler}-\cite{adler-van moerbeke 2}. These papers clearly established  --from a group-theoretical setup-- why standard orthogonality of polynomials and integrability of nonlinear equations of Toda type where so close. In fact, the $LU$ factorization of the moment matrix may be understood as the Gauss--Borel factorization of the initial condition for the integrable hierarchy.  In the Madrid group, based on the Gauss--Borel factorization, we have been searching further the deep links between the Theory of Orthogonal Polynomials and the Theory of Integrable Systems. In \cite{cum1} we studied the generalized orthogonal polynomials \cite{adler} and its matrix extensions from the Gauss--Borel view point. In \cite{cum2} we gave a complete study in terms of factorization for multiple orthogonal polynomials of mixed type and characterized the integrable systems associated to them. Then, we studied Laurent orthogonal polynomials in the unit circle trough the CMV approach in \cite{carlos} and find in \cite{carlos2} the Christoffel--Darboux formula for generalized orthogonal matrix polynomials. These methods where further extended, for example we gave an alternative Christoffel--Darboux formula for mixed multiple orthogonal polynomials \cite{gerardo1} or developed the corresponding  theory of matrix  Laurent orthogonal polynomials in the unit circle and its associated Toda type hierarchy \cite{MOPUC}.

\subsection{Preliminary material}
Here we remind the reader some necessary content extracted from \cite{GM1}. 
Our method to construct Darboux transformations of multivariate  orthogonal polynomials  in a  $D$-dimensional real space (MVOPR) is formulated in terms of a Cholesky factorization of a semi-infinite moment matrix. We consider $D$ independent real variables $\x=\left(x_1,x_2,\dots,x_D \right)^\top\in \Omega\subseteq\mathbb{R}^D$ varying in the domain $\Omega$ together with a Borel measure $\d\mu(x) \in \mathcal{B}(\Omega)$.
The inner product  of two real valued functions $f(\x)$ and $g(\x)$ is defined by
\begin{align*}
\langle f,g\rangle&\coloneq\int_{\Omega} f(\x) g(\x)\d\mu(\x).
\end{align*}

Given a multi-index $\q=(\alpha_1,\dots,\alpha_D)^\top \in\Z_+^{D}$ of non-negative integers we write $\x^{\q}=x_1^{\alpha_1}\cdots x_D^{\alpha_D}$; the length of $\q$ is $|\q|\coloneq  \sum_{a=1}^{D} \alpha_a$. This length induces the total ordering of monomials, $\x^{\q}<\x^{\q'}\Leftrightarrow|\q|<|\q'|$, that we will use to arrange the monomials. For each non-negative integer $k\in\Z_+$ we introduce the set
\begin{align*}
[k]\coloneq \{\q\in \Z_+^{D}: |\q|=k\},
\end{align*}
built up with those vectors  in the lattice $\Z_+^D$ with a given length $k$.
We will use the graded  lexicographic order; i.e.,  for $\q_1,\q_2\in [k]$
\begin{align*}
\q_1>\q_2 \Leftrightarrow \exists p\in \Z_+ \text{ with } p<D \text{ such that } \alpha_{1,1}=\alpha_{2,1},\dots,\alpha_{1,p}=\alpha_{2,p} \text{ and } \alpha_{1,p+1}<\alpha_{2,p+1},
\end{align*}
and if $\q^{(k)}\in[k]$ and $\q^{(\ell)}\in[\ell]$, with $k<\ell$ then $\q^{(k)}<\q^{(\ell)}$.
Given the set of integer vectors of length $k$ we use the  lexicographic order and write
\begin{align*}
[k]=\big\{\q_1^{(k)},\q_2^{(k)},\dots,\q^{(k)}_{|[k]|}\big\} \text{ with } \q_a^{(k)}>\q_{a+1}^{(k)}.
\end{align*}
Here $|[k]|$ is the cardinality of the set $[k]$, i.e., the number of elements in the set.
This is the dimension of the linear space of homogenous multivariate polynomials of total degree $k$.
 Either counting weak compositions or multisets one obtains $|[k]|= \big(\!{D\choose k}\!\big) = {D+k-1 \choose k} $.
The dimension of the linear space $\R_k[x_1,\dots,x_D]$ of multivariate polynomials of degree less or equal  to $k$ is
\begin{align*}
N_{k}=1+|[2]|+\dots+|[k]|=\binom{D+k}{D}.
\end{align*}

We introduce the  vector of monomials
\begin{align*}
\chi&\coloneq \PARENS{\begin{matrix}\chi_{[0]} \\ \chi_{[1]} \\ \vdots \\ \chi_{[k]} \\ \vdots \end{matrix}}
& \mbox{where} & &
\chi_{[k]}&\coloneq  \PARENS{\begin{matrix} \x^{\q_1} \\  \x^{\q_2} \\\vdots \\ \x^{\q_{|[k]|}} \end{matrix}}.
\end{align*}
Observe that for $k=1$ we have that the vectors $\q^{(1)}_a=\ee_a$ for $a\in\{1,\dots,D\}$ forms  the canonical basis of $\R^D$, and for any $\q_j\in[k]$ we have $\q_j=\sum_{a=1}^D \alpha_{j}^a\ee_a$ .
For the sake of simplicity unless needed we will drop off the super-index and write $\q_j$ instead of $\q^{(k)}_j$, as it is understood that $|\q_j|=k$.

The dual space of the symmetric tensor powers is isomorphic to the set of symmetric multilinear functionals on $\R^D$, $\big(\text{Sym}^k(\R^D)\big)^*\cong S((\R^D)^k,\R)$. Hence,
homogeneous polynomials of a given total degree  can be identified with symmetric tensor powers.
Each multi-index $\q\in[k]$ can be thought as a weak $D$-composition of $k$ (or weak composition in  $D$ parts), $k=\alpha_{1}+\dots+\alpha_{D}$.
Notice that these weak compositions may be considered as multisets and that, given a linear basis $\{\ee_a\}_{a=1}^D$ of $\R^D$ we have the linear basis $\{\ee_{a_1}\odot\cdots\odot \ee_{a_k}\}_{\substack{1\leq a_1\leq\cdots\leq a_k\leq D\\ k\in\Z_+}}$ for the symmetric power $\operatorname{S}^k(\R^D)$, where we are using multisets $1\leq a_1\leq\cdots\leq a_k\leq D$. In particular the vectors of this basis $\ee_{a_1}^{\odot M(a_1)}\odot\cdots\odot \ee_{a_p}^{\odot M(a_p)}$, or better its duals $(\ee_{a_1}^*)^{\odot M(a_1)}\odot\cdots\odot (\ee_{a_p}^*)^{\odot M(a_p)}$ are in bijection with monomials of the form $x_{a_1}^{M(a_1)}\cdots x_{a_p} ^{M(a_p)}$. 
The  lexicographic order can be applied  to $\big(\R^D\big)^{\odot k}\cong \R^{|[k]|}$, we then  take a linear basis of $\operatorname{S}^k(\R^D)$ as the ordered set $B_c=\{\ee^{\q_1},\dots,\ee^{\q_{|[k]|}}\}$ with $\ee^{\q_j}\coloneq \ee_1^{\odot \alpha_{j}^1}\odot\dots\odot \ee_{D}^{\odot \alpha_{j}^D}$ so that
$\chi_{[k]}(\x)=\sum_{i=1}^{|[k]|}\x^{\q_j}\ee^{\q_j}$.  

We  consider semi-infinite matrices $A$ with a block or partitioned structure induced by the graded reversed lexicographic order
\begin{align*}
A&=\PARENS{\begin{matrix}
	A_{[0],[0]} & A_{[0],[1]} &  \cdots  \\
	A_{[1],[0]} & A_{[1],[1]} &  \cdots \\
	\vdots                &                 \vdots         &  \\
	\end{matrix}}, &
A_{[k],[\ell]}&=\PARENS{\begin{matrix}
	A_{\q^{(k)}_1,\q^{(\ell)}_1} &   \dots & A_{\q^{(k)}_1,\q^{(\ell)}_{|[\ell]|} }\\
	\vdots & & \vdots\\
	A_{\q^{(k)}_{|[k]|},\q^{(\ell)}_1} &  \dots & A_{\q^{(k)}_{|[k]|},\q^{(\ell)}_{|[\ell]|} }
	\end{matrix}} \in\R^{|[k]|\times |[\ell]|}.
\end{align*}
We use the notation $0_{[k],[\ell]}\in\R^{|[k]|\times|[\ell]|}$ for the rectangular zero matrix, $0_{[k]}\in\R^{|[k]|}$ for the zero vector, and $\I_{[k]}\in\R^{|[k]|\times|[k]|}$ for the identity matrix. For the sake of simplicity we normally  just write $0$ or $\I$ for the zero or identity matrices, and we implicitly assume that the sizes of these matrices are the ones indicated by its position in the partitioned matrix.
\begin{definition}\label{moment}
	Associated with the measure $\d\mu$ we have the following moment matrix
	\begin{align*}
	G&\coloneq \int_{\Omega} \chi(\x)\d\mu(\x) \chi(\x)^\top.
	\end{align*}
	We write the moment matrix  in block form
	\begin{align*}
	G=  \PARENS{\begin{matrix}
		G_{[0],[0]} & G_{[0],[1]} &  \dots \\
		G_{[1],[0]} & G_{[1],[1]} &  \dots \\
		\vdots                &   \vdots              &
		\end{matrix}}.
	\end{align*}
	Truncated  moment matrices are given by
	\begin{align*}
	G^{[\ell]}&\coloneq
	\PARENS{\begin{matrix}
		G_{[0],[0]} &  \cdots & G_{[0],[\ell-1]} \\
		\vdots                        &   & \vdots \\
		G_{[\ell-1],[0]}  &  \cdots & G_{[\ell-1],[\ell-1]}
		\end{matrix}}.
	\end{align*}
\end{definition}
Notice that from the above definition we know that the moment matrix is a symmetric matrix, $G=G^\top$,
which implies that a Gauss--Borel factorization of it, in terms of lower unitriangular \footnote{Lower triangular with the block diagonal populated by identity matrices.} and upper triangular matrices, is a Cholesky factorization.
\begin{pro}\label{qd1}
	If the last quasi-determinants $ \Theta_*(G^{[k+1]})$, $k\in\{0,1,\dots\}$, of the truncated moment matrices are invertible
	the Cholesky factorization
		\begin{align}\label{cholesky}
		G&=S^{-1} H \left(S^{-1}\right)^{\top},
		\end{align}
with
		\begin{align*}
		S^{-1}&=\PARENS{\begin{matrix}
			\I    &             0                &  0                      &  \cdots            \\
			(S^{-1})_{[1],[0]}        & \I&     0                      &   \cdots         \\
			(S^{-1})_{[2],[0]}        & (S^{-1})_{[2],[1]} & \I&      \\
			\vdots                       &        \vdots                  &                               &\ddots
			\end{matrix}}, &
		H&=\PARENS{\begin{matrix}
			H_{[0]}           &   0         &     0         \\
			0                 & H_{[1]} &   0             &    \cdots       \\
			0                  &    0            & H_{[2]} &                        \\
			\vdots   &   \vdots &              &     \ddots       \\
			\end{matrix}},
		\end{align*}
	 can be performed. Moreover,  the rectangular blocks can be expressed in terms of  last quasi-determinants of truncations of the moment matrix
	\begin{align*}
	H_{[k]}&=\Theta_*(G^{[k+1]}), &
	(S^{-1})_{[k],[\ell]}&=\Theta_*(G^{[\ell+1]}_k)\Theta_*(G^{[\ell+1]})^{-1}.
	\end{align*}
\end{pro}
We are ready to introduce the MVOPR
\begin{definition}
	The MVOPR associated to the measure $\d \mu$  are
	\begin{align}\label{eq:polynomials}
	P&=S\chi =\PARENS{\begin{matrix}
		P_{[0]}\\
		P_{[1]}\\
		\vdots
		\end{matrix}}, & P_{[k]}(\x)&=\sum_{\ell=0}^k S_{[k],[\ell]} \chi_{[\ell]}(\x) =\PARENS{\begin{matrix}
		P_{\q^{(k)}_1}\\
		\vdots\\
		P_{\q^{(k)}_{|[k]|}}
		\end{matrix}},&
	P_{\q^{(k)}_i}&=\sum_{\ell=0}^k\sum_{j=1}^{|[\ell]|} S_{\q^{(k)}_i,\q^{(\ell)}_j} \x^{\q^{(\ell)}_j}.
	\end{align}
\end{definition}

Observe that $P_{[k]}=\chi_{[k]}(\x)+\beta_{[k]}\chi_{[k-1]}(\x)+\cdots$ is a vector constructed with the polynomials $P_{\q_i}(\x)$ of degree  $k$, each of which has only one monomial of degree $k$; i. e., we can write $P_{\q_i}(\x)=\x^{\q_i}+Q_{\q_i}(\x)$, with $\deg Q_{\q_i}<k$.
\begin{pro}
	The MVOPR  satisfy
	\begin{align}
	\int_{\Omega} P_{[k]}(\x) \d \mu(\x) (P_{[\ell]}(\x))^{\top}&=\int_{\Omega} P_{[k]} (\x)\d \mu(\x) (\chi_{[\ell]}(\x))^{\top}=0, &
	\ell&=0,1,\dots,k-1,\label{orth}\\
	\int_{\Omega} P_{[k]} (\x)\d \mu(\x)(P_{[k]}(\x))^{\top}&=\int_{\Omega} P_{[k]}(\x) \d\mu(\x)(\chi_{[k]}(\x))^{\top}=H_{[k]}.\label{H}
	\end{align}
\end{pro}
Therefore, we have  the following orthogonality conditions
\begin{align*}
\int_{\Omega} P_{\q^{(k)}_i} (\x)P_{\q^{(\ell)}_j}(\x)\d \mu(\x)&=\int_{\Omega} P_{\q^{(k)}_i}  (\x)\x^{\q^{(\ell)}_j}\d \mu(\x)=0,
\end{align*}
for $\ell=0,1,\dots,k-1$, $i=1,\dots,|[k]|$ and $j=1,\dots,|[\ell]|$,
with the normalization conditions
\begin{align*}
\int_{\Omega} P_{\q_i} (\x)P_{\q_j}(\x)\d \mu(\x)&=\int_{\Omega} P_{\q_i}  (\x)\x^{\q_j}\d \mu(\x)=H_{\q_i,\q_j}, &  i,j&=1,\dots,|[k]|.
\end{align*}
\begin{definition}
	The shift matrices are  given by
	\begin{align*}
	\Lambda_a&=\PARENS{\begin{matrix}
		0   & (\Lambda_a)_{[0],[1]} & 0 & 0 &\cdots\\
		0        & 0 &(\Lambda_a)_{[1],[2]} & 0 &\cdots\\
		0                &  0    &          0       & (\Lambda_a)_{[2],[3]}  &  \\
		0                &  0    &          0             &               0         &\ddots  \\
		\vdots        &  \vdots    &  \vdots         &\vdots
		\end{matrix}}
	\end{align*}
	where the entries in the non zero blocks are given by
	\begin{align*}
	(\Lambda_a)_{\q^{(k)}_i,\q^{(k+1)}_j}&=\delta_{\q^{(k)}_i+\ee_a,\q^{(k+1)}_j},&
	a&=1,\dots, D, &
	i&=1,\dots,|[k]|,&
	j&=1,\dots,|[k+1]|,
	\end{align*}
	and the associated vector 
	\begin{align*}
	 \boldsymbol \Lambda&\coloneq (\Lambda_1,\dots,\Lambda_D)^\top.
	 \end{align*}
	 Finally, we introduce the Jacobi matrices
	 \begin{align}
	 \label{eq:jacobi}
J_a\coloneq &S\Lambda_a S^{-1}, & a\in\{1,\dots,D\},
	 \end{align}
	 and the  vector 
	 \begin{align*}
	 \boldsymbol J=(J_1,\dots,J_D)^\top.
	 \end{align*}
\end{definition}
\begin{pro}\label{pro:Lambda}
	\begin{enumerate}
		\item  The shift matrices commute among them
		\begin{align*}
		\Lambda_a\Lambda_b=\Lambda_b\Lambda_a.
		\end{align*}
		\item We also have the spectral  properties
		\begin{align}\label{eigen}
		\Lambda_a\chi(\x)&= x_a \chi(\x).
		\end{align}
		\item The moment matrix $G$ satisfies
		\begin{align}\label{eq:symmetry}
		\Lambda_a G&= G \big(\Lambda_a\big)^\top.
		\end{align}
		\item The Jacobi matrices $J_a$ are block tridiagonal and satisfy
		\begin{align*}
		J_aH=&HJ_a^\top, & a\in\{1,\dots,D\}.
		\end{align*}
	\end{enumerate}
\end{pro}

Using these properties one derives three term relations or Christoffel--Darboux formul\ae, but as this is not the subject of this paper we refer the interested reader to our paper \cite{GM1}.

\section{Extending the Christoffel formula to the multivariate realm}\label{S:2}
In this section a Darboux transformation for MVOPR is found. Here we use polynomial perturbation of the measure, but to ensure that the procedure works we need perturbations that factor out as $N$ different prime polynomials. Latter we will discuss how we can modify this to include the  most general polynomial perturbation. 
\begin{definition}
Given a degree $m$ polynomial $\mathcal{Q}\in\R[\x]$, $\deg \mathcal Q=m$,
the corresponding  Darboux transformation of the measure is the following perturbed measure
\begin{align*}
\d\mu(\x)\mapsto T\d\mu(\x)\coloneq\mathcal Q(\x)\d\mu(\x).
\end{align*}
\end{definition}
Observe that, if we want a positive definite perturbed measure, we must request to $\mathcal Q$ to be positive definite in the support of the original measure. 

From hereon we assume that both measures $\d\mu$ and $\mathcal Q\d\mu$ give rise to well defined families of MVOPR (equivalently that all their moment matrix block minors are nonzero).

\begin{pro}\label{pro:fac_TG}
	If $TG$ is the moment matrix of $T\d\mu$ we have
	\begin{align*}
	TG=\mathcal Q(\boldsymbol \Lambda)G=G\mathcal (Q(\boldsymbol \Lambda))^\top
	\end{align*}
\end{pro}
\begin{proof}
	It is a direct consequence of the spectral property $\mathcal Q(\boldsymbol{\Lambda})\chi(\x)=\mathcal Q(\x)\chi(\x)$. 
\end{proof}

\begin{definition}
The resolvent matrix is
	\begin{align*}
	\omega\coloneq&(TS)\mathcal Q(\boldsymbol\Lambda)S^{-1}
	\end{align*}
	given in terms of the lower unitriangular matrices $S$ and $TS$ of the Cholesky factorizations of the moment matrices $G=S^{-1}H(S^{-1})^\top$ and  $TG=(TS)^{-1}(TH)(TS^{-1})^\top$.
	The adjoint resolvent is defined by
	\begin{align*}
	M\coloneq S(TS)^{-1}.
	\end{align*}
\end{definition}

\begin{pro}
	We have that
	\begin{align}\label{eq:M-omega}
	(TH)M^\top H^{-1}=\omega.
	\end{align}
\end{pro}
\begin{proof}
	It follows from the Cholesky factorization of $G$ and $TG$ and from \eqref{pro:fac_TG}.
\end{proof}
\begin{pro}\label{pro:superdiagonals}
	In terms of block superdiagonals the resolvent $\omega$ can be expressed as follows
\begin{align*}
  \omega=&\underbracket{\mathcal Q^{(m)}(\boldsymbol{\Lambda})}_{\text{$m$-th superdiagonal}}\\&+
\underbracket{(T\beta)\mathcal Q^{(m-1)}(\boldsymbol{\Lambda})-
\mathcal Q^{(m-1)}(\boldsymbol{\Lambda})\beta}_{\text{$(m-1)$-th superdiagonal}}\\&\shortvdotswithin{+}&+
  \underbracket{(TH)H^{-1}}_{\text{diagonal}}
\end{align*}
\end{pro}
\begin{proof}
The adjoint resolvent is a block lower unitriangular and the resolvent $\omega$ has all its superdiagonals but for the first $m$ equal to zero. The result follows from \eqref{eq:M-omega}.
\end{proof}

\begin{pro}\label{pro:jacobi-LU}
The following $LU$ and $UL$ factorizations
	\begin{align*}
\mathcal Q(\boldsymbol J)=&M\omega,  &
\mathcal Q(T\boldsymbol J)=&\omega M,
\end{align*}
hold.
\end{pro}
\begin{proof}
Both  follow from Proposition \ref{pro:fac_TG} and the Cholesky factorization which imply
\begin{align*}
(TS)^{-1}H(TS^{-1})^\top =\mathcal Q(\boldsymbol \Lambda) S^{-1} H (S^{-1})^\top,
\end{align*}
and a proper cleaning do the job.
\end{proof}
From the first equation in the previous Proposition we get 
\begin{pro}
	The block truncations $(\mathcal Q (\boldsymbol J))^{[k]}$ admit a $LU$ factorization
	\begin{align*}
(\mathcal Q (\boldsymbol J))^{[k]}=M^{[k]}\omega^{[k]}
	\end{align*}
	in terms of the corresponding truncations of the adjoint resolvent $M^{[k]}$ and resolvent $\omega^{[k]}$.
\end{pro}
\begin{pro}\label{pro:regularity-truncation-jacobi}
	We have
	\begin{align*}
	\det (\mathcal Q (\boldsymbol J))^{[k]}=\prod_{l=0}^{k-1}\frac{\det TH_{[l]}}{\det H_{[l]}}
	\end{align*}
and therefore $(\mathcal Q (\boldsymbol J))^{[k]}$ is a regular matrix.
\end{pro}
\begin{proof}
To prove this result just use Propositions \ref{pro:jacobi-LU} and \ref{pro:superdiagonals} and the assumption that
the minors of the moment matrix and the perturbed moment matrix are not zero.
\end{proof}

\begin{pro}\label{pro:the seed}
	The MVOPR  satisfy $\mathcal Q(\x)TP(\x)=\omega P(\x)$. 
	Consequently,
 for any element  $\boldsymbol p$ in the algebraic hypersurface $Z(\mathcal Q)\coloneq \{\x\in\R^D: \mathcal Q(\x)=0\}$ we have
 the important relation
 \begin{align}\label{eq:resolvent}
 \omega_{[k],[k+m]}P_{[k+m]}(\boldsymbol p)+
 \omega_{[k],[k+m-1]}P_{[k+m-1]}(\boldsymbol p)+\cdots+
 \omega_{[k],[k]}P_{[k]}(\boldsymbol p)=0.
 \end{align}
 \end{pro}
 \begin{proof}
 	We have\begin{align*}
 	\omega P(\x)&=(TS)\mathcal Q(\boldsymbol\Lambda)S^{-1}S\chi(\x)\\
 	&=(TS)\mathcal Q(\boldsymbol\Lambda)\chi(\x)\\
 	&=\mathcal Q(\x)(TS)\chi(\x)\\
 	&=\mathcal Q(\x)(TP)(\x).
 	\end{align*}
 
 	Finally, when this formula is evaluated at a point in the algebraic hypersurface of $\mathcal Q$ we obtain that the MVOPR at such points are vectors in the kernel of the resolvent.
 \end{proof}
 To deal with this equation we consider
\begin{definition}
A  set of nodes
\begin{align*}
\mathcal N_{k,m}\coloneq\{\boldsymbol p_j\}_{j=1}^{r_{k,m}}\subset\R^D
\end{align*}
is a set with $r_{k,m}=N_{k+m-1}-N_{k-1}=|[k]|+\dots+|[k+
m-1]|$ vectors in $\R^D$. Given these nodes we consider
the corresponding  sample matrices
\begin{align*}
\Sigma_k^m\coloneq &
\PARENS{
	\begin{matrix}
	P_{[k]}(\boldsymbol p_1) & \dots & P_{[k]}(\boldsymbol p_{r_{k,m}}) \\\vdots& &\vdots\\
	P_{[k+m-1]}(\boldsymbol p_1) & \dots & P_{[k+m-1]}(\boldsymbol p_{r_{k,m}})
	\end{matrix}
}\in\R^{r_{k,m}\times r_{k,m}},\\
\Sigma_{[k,m]}\coloneq&\big(P_{[k+m]}(\boldsymbol p_1), \dots , P_{[k+m]}(\boldsymbol p_{r_{k,m}} )\big)\in\R^{|[k+m]|\times r_{k,m}}.
\end{align*}
\end{definition}

\begin{lemma}\label{lemma:resolvent}
When the set of nodes $\mathcal N_{k,m}\subset Z(\mathcal Q)$ belongs to the algebraic hypersurface of the polynomial $\mathcal Q$ the resolvent coefficients satisfy
\begin{align*}
\omega_{[k],[k+m]}\Sigma_{[k,m]}+
(\omega_{[k],[k]},\dots,\omega_{[k],[k+m-1]})\Sigma_k^m=0.
\end{align*}
\end{lemma}
\begin{proof}
Is a direct consequence of \eqref{eq:resolvent}.
\end{proof}

\begin{definition}
We say that $\mathcal N_{k,m}$  is a poised set if the sample matrix is non singular
\begin{align*}
\det \Sigma_k^m\neq 0.
\end{align*}
\end{definition}

 \begin{theorem}\label{theorem:the deal}
  For a poised set of nodes $\mathcal N_{k,m}\subset Z(\mathcal Q)$ in the algebraic hypersurface of the generating polynomial  $\mathcal Q$  the Darboux transformation of the orthogonal  polynomials can be expressed in terms of the original ones as the following last quasi-determinantal expression
  \begin{align*}
  TP_{[k]}(\x)=  \frac{(\mathcal Q(\boldsymbol\Lambda))_{[k],[k+m]}}{\mathcal Q(\x)}
  \Theta_*\PARENS{\begin{array}{c|c}
  	\Sigma^m_k &\begin{matrix}
  	P_{[k]}(\x) \\ \vdots \\ P_{[k+m-1]}(\x)
  	\end{matrix}\\\hline
  	\Sigma_{[k,m]} & P_{[k+m]}(\x)
  	\end{array}}.
   \end{align*}
 \end{theorem}
 \begin{proof}
Observe that Lemma  \ref{lemma:resolvent} together with $\omega_{[k],[k+m]}=\big(\mathcal Q(\boldsymbol \Lambda)\big)_{[k],[k+m]}$ implies
 \begin{align*}
 (\omega_{[k],[k]},\dots,\omega_{[k],[k+m-1]})=-\big(\mathcal Q(\boldsymbol \Lambda)\big)_{[k],[k+m]}\Sigma_{[k,m]}\big(\Sigma_k^m\big)^{-1}.
 \end{align*}
 and from $\mathcal Q(\x)TP(\x)=\omega P(\x)$ the result follows.
 \end{proof}

\section{Poised sets}\label{S:3}
To construct Darboux transformations in the multivariate setting we need of poised sets in order to find invertible sample matrices with the original polynomials as interpolating functions. When is this possible? Let us start a discussion on this question. First, we introduce two important matrices in the study of poised sets 
  \begin{definition}
We  consider the Vandermonde type matrix
\begin{align*}
\mathcal V_{k}^m\coloneq \big(\chi^{[k+m]}(\boldsymbol p_1) ,\dots,
	\chi^{[k+m]}(\boldsymbol p_{r_{k,m}})\big)\in \R^{N_{k+m}\times r_{k,m}},
\end{align*}
made up of truncated  of multivariate monomials $\chi^{[k+m]}(\x)$ evaluated at the nodes.
We also consider the following truncation  $S_k^m\in\R^{r_{k,m}\times N_{k+m-1}}$ of the lower unitriangular factor $S$ of the Gauss--Borel factorization of the moment matrix
\begin{align}\label{eq:S-slice}
S_k^m\coloneq
\PARENS{
	\begin{matrix}
	S_{[k],[0]} & S_{[k],[1]} &\dots &\I_{|[k]|} & 0_{[k],[k+1]}&\dots & 0_{[k],[k+m-1]}\\
	S_{[k+1],[0]} & S_{[k+1],[1]} &\dots & S_{[k+1],[k]} &\I_{|[k+1]|}&\dots & 0_{[k+1],[k+m-1]}\\
	\vdots & \vdots &  & & & \ddots &\vdots\\
	S_{[k+m-1],[0]} & S_{[k+m-1],[1]} &\dots & & & S_{[k+m-1],[k+m-2]} &\I_{|[k+m-1]|}
	\end{matrix}}.
\end{align}
  \end{definition}
They are relevant because
\begin{lemma}\label{lemma}
We have the following factorization
	\begin{align*}
	\Sigma_k^m=S_k^m
	\mathcal V_{k}^m.
	\end{align*}
\end{lemma}
From where it immediately follows that
\begin{pro}\label{kers}
	The following relations between  linear subspaces
	\begin{align*}
	\operatorname{Ker} \mathcal V^m_k &\subset	\operatorname{Ker} \Sigma_k^m, &
	\operatorname{Im}  \Sigma^m_k &\subset	\operatorname{Im} S_k^m=\R^{r_{k,m}},&
	\end{align*}
	hold true.
\end{pro}

The poisedness of $\mathcal N_{k,m}$ can be reworded as
\begin{align*}
\operatorname{Ker} \Sigma_k^m =\{0\},
\end{align*}
or equivantlently
\begin{align*}
\dim\operatorname{Im}\Sigma_k^m=r_{k,m}.
\end{align*}

\begin{pro}
For poised set $\mathcal N_{k,m}$ the multivariate Vandermonde $\mathcal V_k^m$ needs to have full column rank; i.e., $	\dim\operatorname{Im}\mathcal V_k^m=r_{k,m}$.
\end{pro}
\begin{proof}
 For a set to be poised we need that $\text{Ker}\Sigma_k^m=\{0\}$, but $\text{Ker} \mathcal V_k^m\subset \text{Ker} \Sigma_k^m$ and consequently, 
 \begin{align*}
 \dim\operatorname{Ker}\mathcal V_k^m&=0,
 \end{align*}
 and, as $\dim\operatorname{Ker}\mathcal V_k^m+\dim\operatorname{Im}\mathcal V_k^m=r_{k,m}$, full column rank of the Vandermonde matrix is needed for a set to be poised.
\end{proof}

%
%

The study  of the orthogonal complement of the rank; i.e, the linear subspace $\big(\operatorname{Im}\mathcal V_k^m\big)^\perp\subset \R^{N_{k+m-1}}$ of vectors orthogonal to the image $\operatorname{Im}\mathcal V_k^m$ where $v\in \big(\operatorname{Im}\mathcal V_k^m\big)^\perp$ if
$v^\top\mathcal V_k^m=0$, gives a better insight on the structure of the rank of the Vandermonde matrix.  As $\operatorname{Im}\mathcal V_k^m\oplus\big(\operatorname{Im}\mathcal V_k^m\big)^\perp=\R^{N_{k+m-1}}$ we have the dimensional formula
\begin{align*}
\dim \big(\operatorname{Im}\mathcal V_k^m\big)^\perp +\dim \big(\operatorname{Im}\mathcal V_k^m\big)=N_{k+m-1}.
\end{align*}

\begin{pro}\label{pro:dimensions}
	The Vandermonde matrix $\mathcal V^m_K$ has full column rank if and only if
	\begin{align*}
	\dim \big(\operatorname{Im}\mathcal V_k^m\big)^\perp = N_{k+m-1}-r_{k,m}=N_{k-1}.
	\end{align*}
\end{pro}

Algebraic geometry will guide us in the search of  poised sets  in the algebraic hypersurface of the $m$-th degree polynomial $\mathcal Q$, $\mathcal N_{k,m}\subset Z(\mathcal Q)$.
We need to abandon the real field $\R$  and work in its  algebraical closure $\C$; i.e.,  we understand $\mathcal Q$ as a complex polynomial with real coefficients and consider its zero set as an algebraic hypersurface  in the $D$-dimensional complex affine space $\C^D$; we also change the notation from $\x\in\R^D$ to $\z\in\C^D$.
\begin{definition}
	For a multivariate  polynomial $V$ of total degree $\deg V<k+m-1$  its the principal ideal is $(V)\coloneq\C\{\z^\q V(\z): \q\in\Z_+^D\}\subset\C[z_1,\dots,z_D]$ and for its intersection with the polynomials of degree less or equal than $k+m-1$ we employ the notation
	\begin{align*}
	(V)_{k+m-1}=(V)\cap \C_{k+m-1}[z_1,\dots,z_D]=\C\{\z^\q V(\z)\}_{0\leq|\q|<k+m-\deg V}.
	\end{align*}
\end{definition}
It happens that the elements in the orthogonal complement of the rank of the Vandermonde matrix are polynomials with zeroes at the nodes
\begin{pro}
	As linear spaces the orthogonal complement of the rank of the Vandermonde matrix $\big(\operatorname{Im}\mathcal V_k^m\big)^\perp$ and the space of polynomials of degree less than $k+m$ and  zeroes at  $\mathcal N_{k,m}$ are isomorphic.
\end{pro}
\begin{proof}
	The linear bijection  is
	\begin{align*}
	v=(v_i)_{i=0}^{N_{k+m-1}-1}\in\big(\operatorname{Im}\mathcal V_k^m\big)^\perp&
	\leftrightarrow
	V(\z)=\sum_{i=0}^{N_{k+m-1}-1}v_i \z^{\q_i}
	\end{align*}
	where $V(\z)$ does have zeroes at $\mathcal N_k^m\subset Z(V)$.
	Now, we observe that a vector $v=(v_i)_{i=0}^{N_{k+m-1}-1}\in\big(\operatorname{Im}\mathcal V_k^m\big)^\perp$ can be identified with the   polynomial $V(\z)=\sum_{i}v_i \z^{\q_i}$ which cancels, as a consequence of $v^\top\mathcal V_k^m=0$, at the nodes.  
	\end{proof}
Thus, given this linear isomorphism, for any polynomial $V$ with $\deg V<k+m$ with zeroes  at $\mathcal N_k^m$ we write $V\in \big(\operatorname{Im}\mathcal V_k^m\big)^\perp$. 
\begin{pro}
Given a polynomial $V\in \big(\operatorname{Im}\mathcal V_k^m\big)^\perp$ then
\begin{align*}
(V)_{k+m-1}\subset \big(\operatorname{Im}\mathcal V_k^m\big)^\perp,
\end{align*}
or equivalently 
\begin{align*}
(V)_{k+m-1}^\perp\supseteq \operatorname{Im}\mathcal V_k^m.
\end{align*}
\end{pro}
\begin{proof}
	Observe that all the polynomials in $(V)_{k+m-1}$ have zeroes at the nodes and have degree less or equal to $k+m-1$; hence, we have for the corresponding vectors
	$\{v_\q\}\subset \big(\operatorname{Im}\mathcal V_k^m\big)^\perp$.
\end{proof}

An important result is
\begin{theorem}\label{teorema:ideal-vandermonde}
Given a polynomial $\mathcal Q$, $\deg \mathcal =m$,  and a set of nodes $\mathcal N_{k,m}$	the corresponding Vandermonde matrix $\mathcal V^m_k$ has full column rank if and only if
\begin{align*}
(\mathcal Q)_{k+m-1}=\big(\operatorname{Im}(\mathcal V^m_k)\big)^\perp.
\end{align*}
\end{theorem}
\begin{proof}
	We know that $\dim(\mathcal Q)_{k+m-1}=N_{k-1}$ and as $\mathcal N_{k,m}\subset \mathcal Q$ we know that
	$(\mathcal Q)_{k+m-1}^\perp\supseteq \operatorname{Im}\mathcal V_k^m$. Thus,  full column rankness of the Vandermonde matrix can be achieved if only if  $(\mathcal Q)_{k+m-1}=\big(\operatorname{Im}(\mathcal V^m_k)\big)^\perp$.
\end{proof}
From the spectral property $\mathcal Q(\boldsymbol{\Lambda})\chi(\z)=\mathcal Q(\z)\chi(\z)$ we deduce
\begin{pro}\label{pro:matrix-ideal}
	The row $(\mathcal Q(\boldsymbol\Lambda))_\q$, $\q\in\Z_+^D$, is the graded lexicographic ordering of the entries in the corresponding polynomial $\z^\q\mathcal Q(\z)$.
\end{pro}
Thus, in some way $\mathcal Q(\boldsymbol\Lambda)$ encodes the same information as the principal ideal of $\mathcal Q$ does. To make this  observation  formal we first consider the

\begin{definition}
The matrix  $(\mathcal Q (\boldsymbol\Lambda))^{[k,m]} \in\C^{N_{k-1}\times r_{k,m}}$is  given by
	{\small\begin{align*}\hspace*{-1.55cm}
		(\mathcal Q (\boldsymbol\Lambda))^{[k,m]} \coloneq
		\ccases{	\PARENS{\small
				\begin{matrix}
				(\mathcal Q (\boldsymbol\Lambda))_{[0],[k]}& \dots & (\mathcal Q (\boldsymbol\Lambda))_{[0],[m]}  & 0_{[0],[m+1]} & \dots & 0_{[0],[k+m-1]}\\
				\vdots  &  &&\ddots &\ddots & \vdots & \\
				(\mathcal Q (\boldsymbol\Lambda))_{[k-2],[k]} & & \dots	& &
				& 0_{[k-2],[k+m-1 ]}\\
				(\mathcal Q (\boldsymbol\Lambda))_{[k-1],[k]} && \dots &	& \dots&(\mathcal Q (\boldsymbol\Lambda))_{[k-1],[k+m-1 ]}
				\end{matrix}
			}, &k\leq m,\\
			\PARENS{
				\begin{matrix}
				0_{[0],[k]}&0_{[0],[k+1]} &\dots & 0_{[0],[k+m-1]}\\
				\vdots      &   \vdots&       & \vdots\\
				0_{[k-m-1],[k]} & 0_{[k-m-1],[k+1]}  &\dots & 0_{[k-m-1],[k+m-1]}\\\hline
				(\mathcal Q (\boldsymbol\Lambda))_{[k-m],[k]} & 0_{[k-m],[k+1]} & \dots & 0_{[k-m],[k+m-1]}\\
				(\mathcal Q (\boldsymbol\Lambda))_{[k-m+1],[k]} & 	(\mathcal Q (\boldsymbol\Lambda))_{[k-m+1],[k+1]} & \ddots & \vdots\\
				\vdots & \vdots &\ddots& 0_{[k-2],[k+m-1 ]} \\
				(\mathcal Q (\boldsymbol\Lambda))_{[k-1],[k]} & (\mathcal Q (\boldsymbol\Lambda))_{[k-1],[k+1]}  &\dots& 	(\mathcal Q (\boldsymbol\Lambda))_{[k-1],[k+m-1 ]}
				\end{matrix}
			}, &k\geq m.
		}
		\end{align*}}
	We collect this matrix and the truncation $(\mathcal Q(\boldsymbol{\Lambda}))^{[k]}$ in the $N_{k-1}\times N_{k+m-1}$ matrix
	\begin{align*}
	(\mathcal Q(\boldsymbol{\Lambda}))_k^m\coloneq ((\mathcal Q(\boldsymbol{\Lambda}))^{[k]}, (\mathcal Q(\boldsymbol{\Lambda}))^{[k,m]}).
	\end{align*}
\end{definition}

\begin{pro}\label{pro:idel-vandermonde2}
	We have the following isomorphism 
	\begin{align*}
	P=\sum_{|\q|<k+m}P_\q\z^\q\in(\mathcal Q)_{k+m-1}\Leftrightarrow 
	(P_{\q_0},\dots, P_{\q_{N_{k+m-1}}})=(a_0,\dots, a_{N_{k-1}})(\mathcal Q(\boldsymbol{\Lambda}))_k^m.
	\end{align*}
	between the truncated ideal $(\mathcal Q)_{k+m-1}$ and the orbit of $\C^{N_{k-1}}$ under the linear morphism $	(\mathcal Q(\boldsymbol{\Lambda}))_k^m$. 
	Here we have ordered the multi-indices $\q$ in $P$ according to the graded lexicographic order $\q_0<\dots<\q_{N_{k+m-1}}$.
\end{pro}
\begin{proof}
	Just recall that $P$ is going to be a linear combination of the polynomials $\z^\q\mathcal Q(\z)$, $|\q|<k$.
\end{proof}

Now we  show that full column rankness of the Vandermonde matrix is not only a necessary but also a sufficient condition
\begin{theorem}\label{theorem:poised-fullrank}
Let $\mathcal Q\in\C[\z]$ be a $m$-th degree polynomial, $\deg \mathcal Q=m$. Then, the set of nodes $\mathcal N_{k,m}\subset Z(\mathcal Q)$ is poised if and only if the Vandermonde matrix $\mathcal V^m_k$ has full column  rank.
\end{theorem}
\begin{proof}
Let as assume the contrary, then the sample matrix  $\Sigma_k^m$ is singular, and we can find a nontrivial linear dependence among its rows $(\Sigma_k^m)_i$, $i\in\{1,\dots,r_{k,m}\}$ of the form
\begin{align*}
\sum_{i=1}^{r_{k,m}}c_i(\Sigma_k^m)_i&=0,
\end{align*}
for some nontrivial scalars $\{c_1,\dots,c_{r_{k,m}}\}$.
But, according to Lemma \ref{lemma} $(\Sigma_k^m)_i=(S^m_k)_i\mathcal V^m_k$, where $(S^m_k)_i$ is the $i$-th row of $S^k_m$ and we can write
\begin{align*}
\Big(\sum_{i=1}^{r_{k,m}}c_i(S^m_k)_i\Big)\mathcal V^m_k=0,
\end{align*}
so that
\begin{align*}
\sum_{i=1}^{r_{k,m}}c_i(S^m_k)_i\in\big(\operatorname{Im}\mathcal V^m_k\big)^\perp,
\end{align*}
and given the column full rankness of the Vandermonde matrix, see Theorem \ref{teorema:ideal-vandermonde},
we can write
\begin{align*}
(c_1,\dots,c_{r_{k,m}})S^m_k\in \big(\operatorname{Im}\mathcal V^m_k\big)^\perp,
\end{align*}
or, following Proposition  \ref{pro:idel-vandermonde2} we get
\begin{align}\label{eq:SQ}
(a_0,\dots, a_{N_{k-1}})(\mathcal Q(\boldsymbol{\Lambda}))_k^m=(c_1,\dots,c_{r_{k,m}})S^m_k,
\end{align}
for some non trivial set of $c$'s.
This was the first part of the proof.  

For the second part, we focus on the relation
\begin{align*}
S\mathcal Q(\boldsymbol \Lambda)=\mathcal Q(\boldsymbol J)S.
\end{align*}
Recalling the $S$ is  lower unitriangular by blocks 
\begin{align*}
S=\PARENS{
	\begin{array}{c | c |c c}
	S^{[k]}  & 0  & 0 &\dots\\
	\hline
	\multicolumn{2}{c|} {S^m_k} & 0 &\cdots\\\hline
	\multicolumn{4}{c}{*}
	\end{array}
}
\end{align*}
and that either $\mathcal Q(\boldsymbol \Lambda)$ or $\mathcal Q(\boldsymbol J)$ are block banded matrices with only $m$ block superdiagonals non zero
\begin{align*}
\mathcal Q(\boldsymbol{\Lambda})&=\PARENS{\begin{array}{cc|c}
	(\mathcal Q((\boldsymbol{\Lambda}))^{[k]} & (\mathcal Q((\boldsymbol{\Lambda}))^{[k,m]} &0\\\hline
	*&*&*
	\end{array}}, & 
\mathcal Q(\boldsymbol{J})&=\PARENS{\begin{array}{cc|c}
	(\mathcal Q((\boldsymbol{J}))^{[k]} & (\mathcal Q((\boldsymbol{J}))^{[k,m]} &0\\\hline
	*&*&*
	\end{array}}, 
\end{align*}
which can be written as follows
\begin{align}\label{eq:laleche}
S^{[k]}\big(\mathcal Q(\boldsymbol \Lambda)\big)^m_k=
\big(\mathcal Q(\boldsymbol J)\big)^{[k,m]} S^m_k+\big ((\mathcal Q(\boldsymbol J))^{[k]}S^{[k]},0\big).
\end{align}
We now assume that \eqref{eq:SQ}  holds and multiply \eqref{eq:laleche} by its left with the nonzero row vector \begin{align*}
(a_0,\dots, a_{N_{k-1}})\big(S^{[k]}\big)^{-1}
\end{align*} 
to get
\begin{multline*}
(a_0,\dots, a_{N_{k-1}})\big(\mathcal Q(\boldsymbol \Lambda)\big)^m_k=
(a_0,\dots, a_{N_{k-1}})\big(S^{[k]}\big)^{-1}\big(\mathcal Q(\boldsymbol J)\big)^{[k,m]}  S^m_k\\+
(a_0,\dots, a_{N_{k-1}})\big(S^{[k]}\big)^{-1}\big ((\mathcal Q(\boldsymbol J))^{[k]}S^{[k]},0\big).
\end{multline*}
The expression \eqref{eq:S-slice}  can be written as $S^m_k=(\hat S^m_k,\tilde S^m_k)$ where first block $\hat S^m_k$ is a $r_{k,m}\times N_{k-1}$ matrix and the second block  $\tilde S^m_k$ is a lower unitriangular $r_{k,m}\times r_{k,m}$ matrix. Hence,  the  system \eqref{eq:SQ} can be splitted as follows
\begin{align*}
(a_0,\dots, a_{N_{k-1}})(S^{[k]})^{-1}(\mathcal Q(\boldsymbol J))^{[k,m]}\hat S^m_k+
(a_0,\dots, a_{N_{k-1}})(S^{[k]})^{-1}(\mathcal Q(\boldsymbol J))^{[k]}\big(S^{[k]}\big)^{-1}&=(c_1,\dots,c_{r_{k,m}})\hat S^m_k,\\
(a_0,\dots, a_{N_{k-1}})(S^{[k]})^{-1}(\mathcal Q(\boldsymbol J))^{[k,m]}\tilde S^m_k&=(c_1,\dots,c_{r_{k,m}})\tilde S^m_k.
\end{align*}
As $\tilde S^m_k$ is unitriangular it is invertible and from the second equation we get
\begin{align*}
(a_0,\dots, a_{N_{k-1}})(S^{[k]})^{-1}(\mathcal Q(\boldsymbol J))^{[k,m]}&=(c_1,\dots,c_{r_{k,m}}).
\end{align*}
Thus,  \eqref{eq:SQ}  holds if and only if 
\begin{align*}
(a_0,\dots, a_{N_{k-1}})\big(S^{[k]}\big)^{-1}(\mathcal Q(\boldsymbol J))^{[k]}S^{[k]}=0.
\end{align*}
or equivalently if and only if $\det (\mathcal Q(\boldsymbol J))^{[k]}=0$.
 Thus, recalling Proposition \ref{pro:regularity-truncation-jacobi} (and  our inital assumption that the  two measures $\d\mu(\x)$ and its perturbation $\mathcal Q\d\mu(\x)$ do have MVOPR) the results follows.
\end{proof}


\begin{theorem}\label{H}
Let $\mathcal Q=\mathcal Q_1\cdots \mathcal Q_N$  be the product of $N$ different irreducible polynomials
with $\deg \mathcal Q_a=m_a$, $a\in\{1,\dots,N\}$, and $\deg \mathcal Q=m=\sum\limits_{a=1}^Nm_a$.
Then, the set $\mathcal N_{k,m}\subset Z(\mathcal Q)=\bigcup\limits_{a=1}^N Z(\mathcal Q_a)$ is  poised if the nodes dot not belong to any further complex  algebraic hypersurface of degree smaller than $ k+m$  and different from $Z(\mathcal Q)$.
\end{theorem}
\begin{proof}
Given a subset $Y\subset \C^D$ we define the corresponding ideal $I(Y)=\{P\in\C[z_1,\dots,z_D]: P(\z)=0\, \forall \z\in Y\}$; then, $I\big(\bigcup\limits_{a=1}^N Y_a\big)=\bigcap\limits_{a=1}^N I(Y_a)$ and therefore
$I(Z(\mathcal Q))=\bigcap\limits_{a=1}^NI(Z({\mathcal Q_a}))$. But, according to  the Hilbert's Nullstellensatz and  the prime character  of each factor $\mathcal Q_a$ (every prime ideal is radical) we can write
\begin{align*}
I(Z(\mathcal Q))=&\sqrt{(\mathcal Q)}=\bigcap_{a=1}^N(\mathcal Q_a)\\
=&(\mathcal Q)
\end{align*}
where $\sqrt{(\mathcal Q)}$ is the radical of the principal ideal of $\mathcal Q$. Thus, we conclude
\begin{align*}
\big(\operatorname{Im}\mathcal V_k^m\big)^\perp \supseteq(\mathcal Q)_{k+m-1}
\end{align*}
and deduce
\begin{align*}
\dim\big(\operatorname{Im}\mathcal V_k^m\big)^\perp \geq N_{k-1}.
\end{align*}
The equality is achieved whenever we can ensure that there is no further algebraic hypersurface of degree less than $k+m$, different from $Z(\mathcal Q)$, to which the nodes also belong; i.e., $\operatorname{Im}(\mathcal V^m_k)=(\mathcal Q_1\cdots \mathcal Q_N)_{k+m-1}^\perp$.
\end{proof}
We now discuss on the distribution of nodes along the different irreducible components of the algebraic hypersurface of $\mathcal Q$. 
\begin{pro}
	In a  poised set, $D>1$,  the number $n_a$ of nodes in the irreducible algebraic hypersurface $Z(\mathcal Q_a)$  fulfill
	\begin{align*}
k+m_a \leq n_a&\leq r_{k+m-m_a,m_a}, &
n_1+\cdots+n_N=r_{k,m}.
	\end{align*}
\end{pro}
\begin{proof}
Assume that a number $M_a$ smaller than $r_a\coloneq k+m_a$, $M_a<r_a$, of nodes lay in the irreducible algebraic hypersurface $Z(\mathcal Q_a)$, i.e., the number of nodes in its complementary algebraic hypersurface  $Z(\mathcal Q_1\cdots\mathcal Q_{a-1}\mathcal Q_{a+1}\cdots\mathcal Q_N) $ is bigger than $r_{k,m}-r_a$. Then, the set of nodes belong to the algebraic hypersurface --different of $Z(\mathcal Q)$-- of degree $m-m_a+M_a<k+m$ of the polynomial
$\mathcal Q_1\cdots\mathcal Q_{a-1}\mathcal Q_{a+1}\cdots\mathcal Q_N\pi_1\cdots\pi_{M_a}$, where $\pi_j$ is a degree one polynomial with a zero at the $j$-th node that belongs to  $Z(\mathcal Q_a)$, where we have taken care that $\pi_1\cdots\pi_{M_a}\not\in(\mathcal Q_a)$, which for $D>1$ can be always be achieved. Therefore, we need $M_a\geq k+m_a$ to avoid this situation and to have a  poised set.

 The maximum rank of the  Vandermonde submatrix built up with the columns corresponding to the evaluation of $\chi$ at the nodes in $Z(\mathcal Q_a)$, recalling that  $\dim(\mathcal Q_a)_{k+m-1}=N_{k+m-m_a}$,  is $N_{k+m-1}-N_{k+m-1-m_a}=r_{k+m-m_a,m_a}$.
\end{proof}

Notice, that we need to put $k+m_i$ nodes at each irreducible component $Z(\mathcal Q_a)$, for $a\in\{1,\dots,N\}$, hence we impose conditions on $Nk+m$ nodes. But, do we have enough nodes? The positive answer for $D>1$ can be deduced as follows.
\begin{pro}
The  bound  $r_{k,m}>Nk+m$ holds.
\end{pro}
\begin{proof}
We have $r_{k,m}=|[k+m-1]|+\cdots+|[k]|$, thus a  rude lower bound of nodes  (for $D>1$) is
\begin{align*}
r_{k,m}>&m|[k]|=m\binom{k+D-1}{D-1}=m\Big(1+\frac{k}{D-1}\Big)\cdots\Big(1+\frac{k}{2}\Big)(1+k)\\
>&m(k+1)\\
>&Nk+m.
\end{align*}
\end{proof}

But, what happens with this condition for $D=1$? Now, we have $r_{k,m}=m$ nodes, $m_a=1$, each $Z(\mathcal Q_i)$ is a single point in $\C$ and $N=m$. In this case, the reasoning that lead to the construction of the polynomial $\mathcal Q_1\cdots\mathcal Q_{a-1}\mathcal Q_{a+1}\cdots\mathcal Q_N\pi_1\cdots\pi_{M_a}$ in the previous proof is not applicable; first $M_a=1$, given that all prime factors are degree one polynomials and second, the polynomial $\pi_1$ must be $\mathcal Q_a$ and therefore the product leads to the polynomial $\mathcal Q$ and no further constraint must be considered.

 The maximum $n_a$ is greater than the minimum number of nodes of that type $r_{k+m-m_a,m_a}>m_a(k+m-m_a+1)>k+m_a$.  Moreover, the sum of the maximum ranks exceeds the number of nodes, and  the full column rank condition is reachable:
\begin{pro}
We have $\sum\limits_{i=1}^{N}r_{k+m-m_a,m_a}\geq r_{k,m}$.
\end{pro}
\begin{proof}
	 For $N=2$ we need to show that
	 $r_{k+m_2,m_1}+r_{k+m_1,m_2}>r_{k,m_1+m_2}$ or
	 \begin{multline*}
	 |[k+m_2+m_1-1]|+\cdots+|[k+m_2]|+|[k+m_2+m_1-1]|+\cdots+|[k+m_1]|\\>
	 |[k+m_2+m_1-1]|+\cdots+|[k+m_2]|+|[k+m_2-1]|+\cdots+|[k]|
	 \end{multline*}
	 which is obvious. Then, for $N=3$ we need to prove that $r_{k+m_2+m_3,m_1}+r_{k+m_1+m_3,m_2}+r_{k+m_1+m_2,m_3}>r_{k,m_1+m_2+m_3}$, but
	 using the  already proven $N=2$ case we have $r_{k+m_2+m_3,m_1}+r_{k+m_1+m_3,m_2}>r_{k+m_3,m_1+m_2}$ and using the $N=2$ equation again we do have $r_{k+m_3,m_1+m_2}+r_{k+m_1+m_2,m_3}>r_{k,m_1+m_2+m_3}$, as desired. An induction procedure gives the result for arbitrary $N$.
\end{proof}

In the next picture we illustrate the case $N=2$ of two prime polynomials of degrees $m_1$ and $m_2$. The blue rectangle gives the possible values for the number of nodes $(n_1,n_2)$ corresponding $n_i$ to the prime polynomial $\mathcal Q_i$ according to the bounds
\begin{align*}
k+m_1 \leq n_1&\leq r_{k+m_2,m_1}, &
k+m_2 \leq n_2&\leq r_{k+m_1,m_2}.
\end{align*}
The blue diagonals $n_1+n_2=K$ are ordered according  (we assume for the degrees  that $m_1\leq m_2$ and therefore $r_{k+m_2,m_1}\leq r_{k+m_1,m_2}$)  to the chain of inequalities
\begin{align*}
2k+m\leq  \max(r_{k+m_2,m_1}+k+m_2,r_{k+m_1,m_2}+k+m_1)\leq r_{k,m}\leq r_{k+m_2,m_1}+r_{k+m_1,m_2}.
\end{align*}
Notice that $\max(r_{k+m_2,m_1}+k+m_2,r_{k+m_1,m_2}+k+m_1)\leq r_{k,m}$ follows $r_{k,m}=r_{k+m_i,m_j}+r_{k,m_i}\geq r_{k+m_i,m_j}+m_i|[k]|\geq  r_{k+m_i,m_j}+m_i(k+1)\geq   r_{k+m_i,m_j}+k+m_i$, where  $(i,j)=(1,2),(2,1)$. Therefore, the striped triangle is the area where the couples $(n_1,n_2)$ of number of nodes belong. We have drawn the passing of the line $n_1+n_2=r_{k,m}$ trough it, and show the integer couples in that segment, those will be the possible distributions of nodes among the zeroes of both prime polynomials.
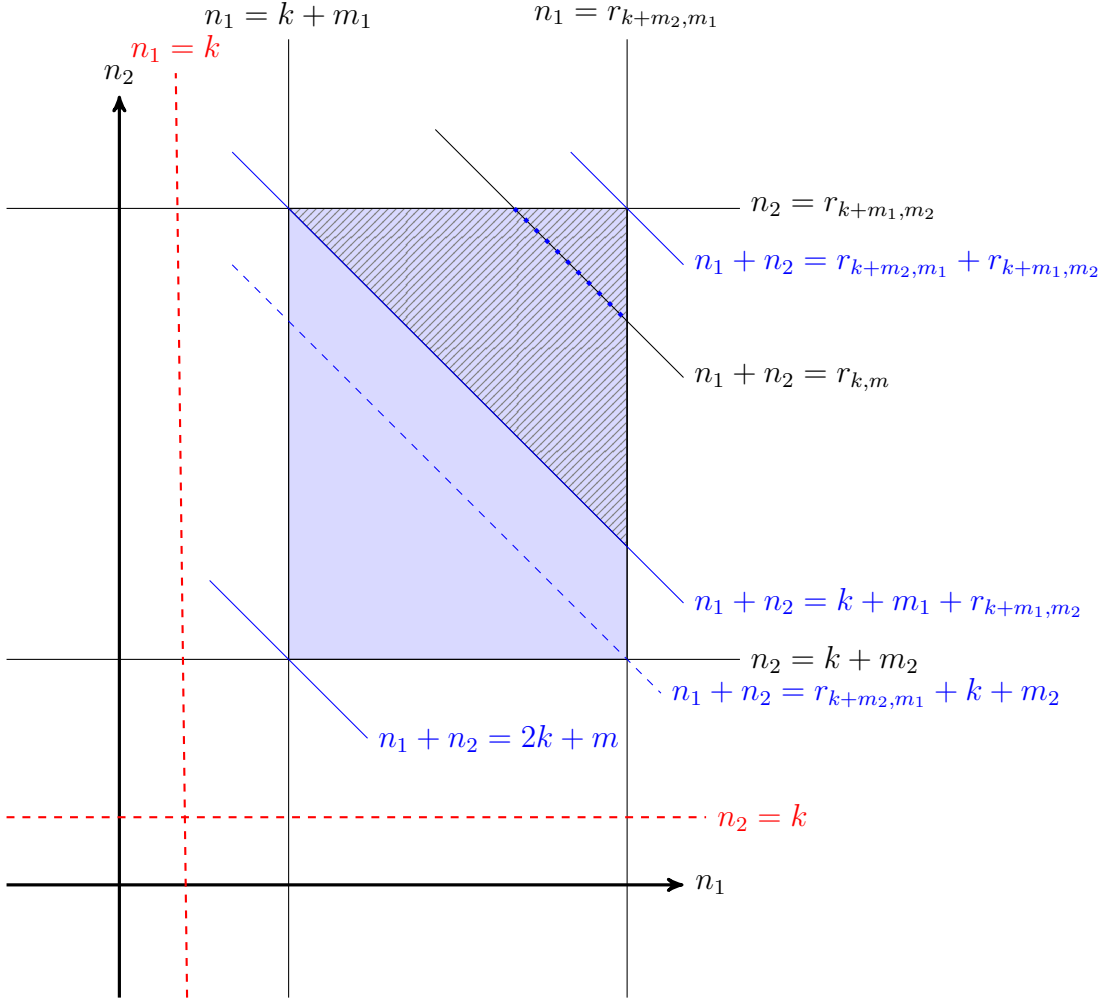
\begin{center}
\begin{tikzpicture}[
scale=1.5,
axis/.style={very thick, ->, >=stealth'},
]
\draw[axis] (-1,0)  -- (5,0) node(xline)[right] {$n_1$};
\draw[axis] (0,-1) -- (0,7) node(yline)[above] {$n_2$};
\filldraw[fill=blue!15!white, draw=black] (1.5,2) rectangle (4.5,6);
\filldraw[pattern=north east lines, pattern color=gray] (4.5,6)--(4.5,3)--(1.5,6);
\draw (1.5,-1) -- (1.5,7.5)  node[above] {$n_1=k+m_1$};
\draw (4.5,-1) -- (4.5,7.5) node[above] {$n_1=r_{k+m_2,m_1}$};
\draw (-1,2) -- (5.5,2) node[right] {$n_2=k+m_2$};
\draw (-1,6) -- (5.5,6) node[right] {$n_2=r_{k+m_1,m_2}$};
\draw[red,thick,dashed]  (0.6,-1) -- (0.5,7.2) node[above] {$n_1=k$} ;
\draw[red,thick,dashed]  (-1,0.6) -- (5.2,0.6) node[right] {$n_2=k$} ;
\draw[color=blue, domain=0.8:2.2]    plot (\x,-\x+3.5)             node[right] {$n_1+n_2=2k+m$};
\draw[color=blue, domain=1:5]    plot (\x,-\x+7.5)    node[right] {$n_1+n_2=k+m_1+r_{k+m_1,m_2}$};
\draw[color=blue, domain=1:4.8,dashed]    plot (\x,-\x+6.5)             node[right] {$n_1+n_2=r_{k+m_2,m_1}+k+m_2$};
\draw[color=blue, domain=4:5]    plot (\x,-\x+10.5)             node[right] {$n_1+n_2=r_{k+m_2,m_1}+r_{k+m_1,m_2}$};
\draw[domain=2.8:5,thin]    plot (\x,-\x+9.5)             node[right] {$n_1+n_2=r_{k,m}$};
\draw[loosely dotted,domain=3.5:4.51,blue,ultra thick]    plot (\x,-\x+9.5);
\end{tikzpicture}
\end{center}

\section{Darboux transformations for a general perturbation}\label{S:4}
We begin with a negative result 
\begin{pro}
	Poised sets do not exist for $\mathcal Q=\mathcal R^d$, $d\in\{2,3,\dots\}$, for any given  polynomial $\mathcal R$.
\end{pro}
\begin{proof}
	Now $\mathcal Q=\mathcal R^d$, $d\in\{2,3,4,\dots\}$  $\deg \mathcal Q=d \deg \mathcal R$, for some polynomial $\mathcal R$. In this case $Z(\mathcal Q)=Z({\mathcal R})$, but $\dim(\mathcal R)_{k+m-1}=N_{k-1+(d-1)\deg\mathcal R}>N_k$ and consequently the set is not  poised.
\end{proof}

We now discuss a  method to overcome this situation. We will generalize the construction of nodes, sample matrices and poised sets. In this manner we are able to give explicit Christoffel type formulae for the Darboux transformation of more general generating polynomials. 
 We consider multi-Wro\'{n}ski type matrices and multivariate confluent  Vandermonde matrices.

\subsection{Discussion for the arbitrary power of a  prime polynomial}
Now we take $\mathcal Q=\mathcal R^d$,  $\deg\mathcal R=n$ and $\deg\mathcal Q=dn$, so that $Z(\mathcal Q)=Z(\mathcal R)$ with $\mathcal R$ to be a prime polynomial.  From Proposition \ref{pro:the seed} we know that
$\omega P(\x)= \mathcal R^d(\x)TP(\x)$.
To analyze this situation we consider  a set of linearly independent vectors $\big\{\n^{(j)}_i\big\}_{i=1}^{\rho_j}\subset \R^{|[j]|}\cong \big(\R^D\big)^{\odot j}$,  $\rho_j\leq |[j]|$; here $\n_i^{(j)}=(n^{(j)}_{i,\q})_{\q\in [j]}$, and to each of these vectors we  associate the following homeogenous linear  differential operator
\begin{align*}
\frac{\partial^j}{\partial \n_i^{(j)}}=\sum_{|\q|=j}n^{(j)}_{i,\q}\frac{\partial^j}{\partial \x^\q}.
\end{align*}

From the Leibniz rule we infer
\begin{pro}\label{resolventd}
	For any element  $\boldsymbol p$ in the algebraic hypersurface $Z(\mathcal R)\coloneq \{\x\in\R^D: \mathcal R(\x)=0\}$ we have
	\begin{align*}
	\omega_{[k],[k+nd]}\frac{\partial^jP_{[k+nd]}}{\partial \n_i^{(j)}}(\boldsymbol p)+
	\omega_{[k],[k+nd-1]}\frac{\partial^jP_{[k+nd-1]}}{\partial \n_i^{(j)}}(\boldsymbol p)+\cdots+
	\omega_{[k],[k]}\frac{\partial^jP_{[k]}}{\partial \n_i^{(j)}}(\boldsymbol p)&=0,
	\end{align*}
	for $j\in\{0,1,\dots,d-1\}$ and $i\in\{1,\dots,\rho_j\}$.
\end{pro}
This suggests to  extend the set of nodes and the  sample matrices 
\begin{definition}
	We consider the splitting into positive integers $r_{k,nd}=N_{k+n d-1}-N_{k-1}=|[k]|+\dots+|[k+
	d n-1]|=\sum\limits_{j=0}^{d-1}\sum\limits_{i=1}^{\rho_j}\nu_i^{(j)}$, and for each $j\in\{0,1,\dots,d-1\}$ we consider  the following set of distinct nodes
	\begin{align*}
	\mathcal N_{i}^{(j)}\coloneq\big\{\boldsymbol p^{(j)}_{i,l}\big\}_{l=1}^{\nu_i^{(j)}}\subset\R^D,
	\end{align*}
	where we allow for non empty intersections  between these sets of nodes and we denote its union by $\mathcal {N}_{k,nd}=\bigcup\limits_{j=0}^{d-1}\bigcup\limits_{i=1}^{\rho_j}\mathcal {N}_i^{(j)}$.
	We also need of the above mentioned  set of linearly independent vectors $\big\{\n^{(j)}_i\big\}_{i=1}^{\rho_j}\subset \R^{|[j]|}\cong \big(\R^D\big)^{\odot j}$,  $\rho_j\leq |[j]|$, $j\in\{0,\dots,d-1\}$.
	The  partial blocks of the homogeneous sample matrices  are
	\begin{align*}
	(\Sigma_k^{nd})^{(j)}_i \coloneq&
		\PARENS{
			\begin{matrix}
			\dfrac{\partial^j P_{[k]}}{\partial\n_i^{(j)}}(\boldsymbol p^{(j)}_{i,1}) & \dots & 	\dfrac{\partial^j P_{[k]}}{\partial\n_i^{(j)}}(\boldsymbol p^{(j)}_{i,\nu_i^{(j)}}) \\\vdots& &\vdots\\
			\dfrac{\partial^j P_{[k+nd-1]}}{\partial\n_i^{(j)}}(\boldsymbol p^{(j)}_{i,1}) & \dots & 	\dfrac{\partial^j P_{[k+nd-1]}}{\partial^j\n_i^{(j)}}(\boldsymbol p^{(j)}_{i,\nu_i^{(j)}})
			\end{matrix}   }\in\R^{r_{k,nd}\times \nu^{(j)}_i},
	\\
	(\Sigma_{[k,nd]})^{(j)}_i\coloneq&
		\Big(
		\dfrac{\partial^j P_{[k+nd]}}{\partial\n_i^{(j)}}(\boldsymbol p^{(j)}_{i,1}),  \dots , 	\dfrac{\partial^j P_{[k+nd]}}{\partial\n_i^{(j)}}(\boldsymbol p^{(j)}_{i,\nu_i^{(j)}})
		\Big)
		\in\R^{|[k+nd]|\times  \nu^{(j)}_i},
	\end{align*}
	in terms of which we write the homogenous sample matrices 
	\begin{align*}
	\big(\Sigma_k^{nd}\big)^{(j)}\coloneq&\big(\big(\Sigma_k^{nd}\big)^{(j)}_1,\dots,\big(\Sigma_k^{nd}\big)^{(j)}_{\rho_j} \big)\in\R^{r_{k,nd}\times \sum_{i=1}^{\rho_j}\nu_i^{(j)}},\\\big(\Sigma_{[k,nd]}\big)^{(j)}\coloneq&
	\big(\big(\Sigma_{[k,nd]}\big)^{(j)}_1,\dots,\big(\Sigma_{[k,nd]}\big)^{(j)}_{\rho_j} \big)\in\R^{|[k+nd]|\times \sum_{i=1}^{\rho_j}\nu_i^{(j)}}, 
	\end{align*}
	which allow us to define  the  multivariate Wro\'{n}ski type sample matrices 
	\begin{align*}
	\Sigma_k^{nd}\coloneq&\big((\Sigma_k^{nd})^{(0)},  ,\dots,(\Sigma_k^{nd})^{(d-1)} \big)\in\R^{r_{k,nd}\times r_{k,nd}},\\
	\Sigma_{[k,nd]}\coloneq&\big((\Sigma_{[k,nd]})^{(0)},\dots,(\Sigma_{[k,nd]})^{(d-1)}\big)\in\R^{|[k+nd]|\times r_{k,nd}}.
	\end{align*}
\end{definition}

\begin{definition}
	We say that $\mathcal {N}_{k,nd}$  is a poised set if the sample matrix is non singular
	\begin{align*}
	\det\Sigma_k^{nd}\neq 0.
	\end{align*}
\end{definition}

\begin{theorem}
	For a poised set of nodes $\mathcal {N}_{k,nd}\subset Z(\mathcal Q)$ in the algebraic hypersurface of the prime polynomial $\mathcal R$  the transformed orthogonal polynomials can be expressed in terms of the original ones as according to the  quasi-determinantal expression
	\begin{align*}
	TP_{[k]}(\x)=  \frac{(\mathcal R(\boldsymbol\Lambda)^d)_{[k],[k+nd]}}{\mathcal R(\x)^d}
	\Theta_*\PARENS{\begin{array}{c|c}
		\Sigma^{nd}_k &\begin{matrix}
		P_{[k]}(\x) \\ \vdots \\ P_{[k+nd-1]}(\x)
		\end{matrix}\\\hline
		\Sigma_{[k,nd]} & P_{[k+nd]}(\x)
		\end{array}}.
	\end{align*}
\end{theorem}
\begin{proof}
	Proposition \ref{resolventd} gives
	\begin{align*}
	\omega_{[k],[k+nd]}\Sigma_{[k,nd]}+
	\big(\omega_{[k],[k]},\cdots,	\omega_{[k],[k+nd-1]}\big)\Sigma_k^{nd}=0
	\end{align*}
so that
	\begin{align*}
	(\omega_{[k],[k]},\dots,\omega_{[k],[k+nd-1]})=-\big(\mathcal R(\boldsymbol \Lambda)^d\big)_{[k],[k+nd]}\Sigma_{[k,nd]}\big(\Sigma_k^{nd}\big)^{-1}.
	\end{align*}
	and  $\mathcal R(\x)^dTP(\x)=\omega P(\x)$ gives the result.
\end{proof}
To discuss the existence of  poised sets we allow the nodes to be complex.
\begin{definition}
	We  introduce the partial derived Vandermonde matrices
	\begin{align*}
	(\mathcal V_{k}^{nd})^{(j)}_i\coloneq&
		\Big(\dfrac{\partial^j \chi^{[k+nd]}}{\partial \n^{(j)}_i}(\boldsymbol p^{(j)}_{i,1}) ,\dots,
		\dfrac{\partial^j \chi^{[k+nd]}}{\partial \n^{(j)}_i}(\boldsymbol p^{(j)}_{i,\nu^{(j)}_i})\Big)\in \C^{N_{k+nd-1}\times \nu^{(j)}_i},
	\end{align*}
	for $ j\in\{0,\dots,d-1\}$ and $i\in\{1,\dots,\rho_j\}$,
	 the derived  Vandermonde matrix is
	 \begin{align*}
	 (\mathcal V_{k}^{nd})^{(j)}\coloneq&
	 \big((\mathcal V_{k}^{nd})^{(j)}_1,\dots,(\mathcal V_{k}^{nd})^{(j)}_{\rho_j}\big)\in\C^{N_{k+nd-1}\times \sum_{i=1}^{\rho_j}\nu^{(j)}_i},
	 \end{align*}
	 and the multivariant confluent Vandermonde matrix
	\begin{align*}
	\mathcal V^{nd}_k\coloneq \big( (\mathcal V_{k}^{nd})^{(0)},(\mathcal V_{k}^{nd})^{(1)}\dots,(\mathcal V_{k}^{nd})^{(d-1)}\big)\in\C^{N_{k+nd-1}\times r_{k,nd}}.
	\end{align*}
\end{definition}
As in the previous analysis we have
	$\Sigma_k^m=S_k^m
	\mathcal V_{k}^m$,
	where $S^m_k$ is given in \eqref{eq:S-slice},
and $\operatorname{Ker}\mathcal V^{nd}_k \subset	\operatorname{Ker} \Sigma_k^{nd}$.
For $\mathcal N_{k,nd}$ to be  poised  we must request to  $\mathcal V_k^{nd}$ to be a full column rank matrix; i.e.,
$\dim\operatorname{Im}\mathcal V_k^{nd}=r_{k,nd}$. Remarkably, Theorems \ref{teorema:ideal-vandermonde} and \ref{theorem:poised-fullrank} and Propositions  \ref{pro:matrix-ideal} and \ref{pro:idel-vandermonde2} hold true for our polynomial $\mathcal Q=\mathcal R^d$ and the corresponding confluent Vandermonde matrix $\mathcal V^m_k$, $m=nd$.


\begin{theorem}
	The node set $\mathcal N_{k,nd}\subset\C^D\subset Z(\mathcal Q)$ is   poised if it does not exist a polynomial  $V\neq \mathcal R^d$, $\deg V\leq k+nd-1$, such that 
	$\mathcal N_j^{(i)}\subset Z\Big(\dfrac{\partial^j V}{\partial \n^{(j)}_i}\Big)$, 	for $ j\in\{0,\dots,d-1\}$ and $i\in\{1,\dots,\rho_j\}$.
\end{theorem}
\begin{proof}
A vector $v=(v_i)_{i=1}^{N_{k+nd-1}}\in(\operatorname{Im}\mathcal V_k^{nd})^\perp$ if for the corresponding polynomial $V=\sum\limits_{l=1}^{N_{k+nd-1}}v_l\x^{\q_l}$  the polynomials $\dfrac{\partial^j V}{\partial\n^{(j)}_i}$ cancel at $\mathcal N^{(j)}_i$. Remarkably,   $\dfrac{\partial^j (x^\q\mathcal R^d)}{\partial\n^{(j)}_i}(\boldsymbol p)=0$, $j=1,\dots,d-1$ for $i\in\{1,\dots,\rho_j\}$ and $\boldsymbol p\in Z(\mathcal R)$. Hence, we conclude that
	\begin{align*}
	(\mathcal R^d)_{k+nd-1}\subseteq (\operatorname{Im}\mathcal V_k^{nd})^\perp,
	\end{align*}
	and, as $\dim (\mathcal R^d)_{k+nd-1}=N_{k-1}$, we have full column rankness of the confluent Vandermonde matrix  if  $\operatorname{Im}\mathcal V^{nd}_k=(\mathcal R^d)_{k+nd-1}^\perp$.
\end{proof}

\begin{cor}For a  poised set
	\begin{itemize}
		\item we can not take the vectors $\n_i^{(j)}\in\R^{|[j]|}$ such that for a given $p\in\{1,\dots,d\}$ the polynomial $\dfrac{\partial^j (\mathcal R^p)}{\partial\n^{(j)}_i}$  cancels  at $Z(\mathcal R)$.
		\item  we can not pick up the  nodes from an algebraic hypersurface of degree less than or equal to $\lfloor\frac{k-1}{d}\rfloor+n$.\footnote{Here $\lfloor x\rfloor$ is the floor function and gives the greatest integer less than or equal to $x$.}
		\item the following upper bounds must hold
		\begin{align*}
		\nu_0^{(0)}&\leq r_{k+n(d-1),n}, & \nu_i^{(j)}&\leq r_{k+nd-d_i^{(j)},d_i^{(j)}},
		\end{align*}
		where $d_i^{(j)}\coloneq\deg\dfrac{\partial(\mathcal R^d)}{\partial \n^{(j)}_i}$.
	\end{itemize}
\end{cor}
\begin{proof}
When $\dfrac{\partial ^j(\mathcal R^p)}{\partial\n_i^{(j)}}$  cancels at $Z(\mathcal R)$ then $\dim (\operatorname{Im}\mathcal V_k^m)^\perp\geq N_{k+(d-p)n}$ and
	the set is not  poised. Given a polynomial $W$, $\deg W\leq \lfloor\frac{k-1}{d}\rfloor+n$, of the  described type we see that $V=W^d$, $\deg V\leq k-1+m$, is a polynomial such that $\frac{\partial^j V}{\partial\n_i^{(j)}}$ cancels at $Z(W)$ and again full column rankness is not achievable.
	All the columns  in the  Vandermonde block $(\mathcal V_k^{nd})^{(0)}$  (we have remove the subindex because for $j=0$ there is only one and no need to distiguish among several of them) imply no directional partial derivatives, so that  $(\mathcal R)_{k+nd-1}^\perp\supseteq\operatorname{Im}(\mathcal V_k^{(nd)})^{(0)}$ and the maximum achievable rank  for this block is $N_{k+nd-1}-N_{k+n(d-1)-1}=r_{k+n(d-1),n}$.
For $j=1,\dots,d-1$ the columns used in the construction of  the block $(\mathcal V_k^{nd})^{(j)}_i$
imply directional partial derivatives $\frac{\partial^j\quad}{\partial \n^{(j)}_i}$
	and consequently  $\big(\frac{\partial^j(\mathcal R^d))}{\partial \n^{(j)}_i}\big)_{k+nd-1}^\perp\supseteq \Big(\operatorname{Im}\big((\mathcal V_k^{nd})^{(j)}_i\big)\Big)$; hence, the maximum rank is $N_{k+nd-1}-N_{k+nd-1-d^{(j)}_i}$.
	\end{proof}

\subsection{The general case}

We now consider the general situation of a polynomial in several variables, i.e., $\mathcal Q=\mathcal R_1^{d_1}\cdots \mathcal R_N^{d_N}$ where $\mathcal R_i$, $\deg \mathcal R_i=m_i$, are different prime polynomials; we have for the degree of the polynomial $\deg\mathcal Q=m=n_1d_1+\cdots+n_Nd_N$. As with the study of the product of $N$ different prime polynomials developed in \S. \ref{S:3}
we have $Z(\mathcal Q)=\bigcup\limits_{i=1}^NZ(\mathcal R_i)$.

\begin{definition}
	 Consider the splitting $r_{k,n_1d_1+\cdots +n_Nd_N}=N_{k+n_1d_1+\cdots +n_Nd_N-1}-N_{k-1}=|[k]|+\dots+|[k+
	dn_1d_1+\cdots +n_Nd_N-1]|=\sum\limits_{a=1}^N\sum\limits_{j=0}^{d_a-1}\sum\limits_{i=1}^{\rho_{a,j}}\nu_i^{(a,j)}$, and for each $j\in\{0,1,\dots,d_a-1\}$ the following set of different nodes
	\begin{align*}
	\mathcal N_{i}^{(a,j)}\coloneq&\big\{\boldsymbol p^{(a,j)}_{i,l}\big\}_{l=1}^{\nu_i^{(a,j)}}\subset\R^D,&
		\mathcal N_a\coloneq&\bigcup\limits_{j=0}^{d_a-1}\bigcup\limits_{i=1}^{\rho_j}\mathcal {N}_i^{(a,j)}
	\end{align*}
	where for a fixed $a\in\{1,\dots,N\}$ we allow for non empty intersections among  sets with different values of $j$; denote its union by $\mathcal {N}_{k,n_1d_1+\cdots +n_Nd_N}=\bigcup\limits_{a=1}^N\mathcal N_a$ and 
pick a  set of linearly independent vectors $\big\{\n^{(a,j)}_i\big\}_{i=1}^{\rho_{a,j}}\subset \R^{|[j]|}\cong \big(\R^D\big)^{\odot j}$,  $\rho_j\leq |[j]|$, $j\in\{0,\dots,d_a-1\}$.
	The associated  $i$-th  homogeneous blocks of the sample matrices are
	\begin{align*}
	(\Sigma_k^{n_1d_1+\cdots +n_Nd_N})^{(a,j)}_i \coloneq&
		\PARENS{
			\begin{matrix}
			\dfrac{\partial^j P_{[k]}}{\partial\n_i^{(a,j)}}(\boldsymbol p^{(a,j)}_{i,1}) & \dots & 	\dfrac{\partial^j P_{[k]}}{\partial\n_i^{(a,j)}}(\boldsymbol p^{(a,j)}_{i,\nu_i^{(a,j)}}) \\\vdots& &\vdots\\
			\dfrac{\partial^j P_{[k+n_1d_1+\cdots +n_Nd_N-1]}}{\partial\n_i^{(a,j)}}(\boldsymbol p^{(a,j)}_{i,1}) & \dots & 	\dfrac{\partial^j P_{[k+n_1d_1+\cdots +n_Nd_N-1]}}{\partial^j\n_i^{(a,j)}}(\boldsymbol p^{(a,j)}_{i,\nu_i^{(a,j)}})
			\end{matrix}   },
	\\
	(\Sigma_{[k,n_1d_1+\cdots +n_Nd_N]})^{(a,j)}_i\coloneq&
		\Big(
		\dfrac{\partial^j P_{[k+n_1d_1+\cdots +n_Nd_N]}}{\partial\n_i^{(a,j)}}(\boldsymbol p^{(a,j)}_{i,1}),  \dots , 	\dfrac{\partial^j P_{[k+n_1d_1+\cdots +n_Nd_N]}}{\partial\n_i^{(a,j)}}(\boldsymbol p^{(a,j)}_{i,\nu_i^{(a,j)}})
		\Big),
	\end{align*}
with $(\Sigma_k^{n_1d_1+\cdots +n_Nd_N})^{(a,j)}_i\in\R^{r_{k,n_1d_1+\cdots +n_Nd_N}\times \nu^{(a,j)}_i}$ and 	$(\Sigma_{[k,n_1d_1+\cdots +n_Nd_N]})^{(a,j)}_i\in\R^{|[k+n_1d_1+\cdots +n_Nd_N]|\times  \nu^{(a,j)}_i}$,
	 the homogenous sample matrices are
	\begin{align*}
	\big(\Sigma_k^{n_1d_1+\cdots +n_Nd_N}\big)^{(a,j)}\coloneq&\big(\big(\Sigma_k^{n_1d_1+\cdots +n_Nd_N}\big)^{(a,j)}_1,\dots,\big(\Sigma_k^{n_1d_1+\cdots +n_Nd_N}\big)^{(a,j)}_{\rho_{a,j}} \big)\in\R^{r_{k,n_1d_1+\cdots +n_Nd_N}\times \sum_{i=1}^{\rho_{a,j}}\nu_i^{(a,j)}},\\\big(\Sigma_{[k,n_1d_1+\cdots +n_Nd_N]}\big)^{(a,j)}\coloneq&
	\big(\big(\Sigma_{[k,n_1d_1+\cdots +n_Nd_N]}\big)^{(a,j)}_1,\dots,\big(\Sigma_{[k,n_1d_1+\cdots +n_Nd_N]}\big)^{(a,j)}_{\rho_{a,j}} \big)\in\R^{|[k+n_1d_1+\cdots +n_Nd_N]|\times \sum_{i=1}^{\rho_{a,j}}\nu_i^{(a,j)}}
	\end{align*}
	and the partial multivariate Wro\'{n}ski type sample matrices are
	\begin{align*}
	(\Sigma_k^{n_1d_1+\cdots +n_Nd_N})_a\coloneq&\big((\Sigma_k^{n_1d_1+\cdots +n_Nd_N})^{(a,0)},  ,\dots,(\Sigma_k^{n_1d_1+\cdots +n_Nd_N})^{(a,d_a-1)} \big),\\
	(\Sigma_{[k,n_1d_1+\cdots +n_Nd_N]})_a\coloneq&\big((\Sigma_{[k,n_1d_1+\cdots +n_Nd_N]})^{(a,0)},\dots,(\Sigma_{[k,n_1d_1+\cdots +n_Nd_N]})^{(a,d-1)}\big),
	\end{align*}
that are matrices in $\R^{r_{k,n_1d_1+\cdots +n_Nd_N}\times \sum_{j=0}^{d_a-1}\sum_{i=1}^{\rho_{a,j}}\nu_i^{(a,j)}}$ and in $\R^{|[k+n_1d_1+\cdots +n_Nd_N]|\times \sum_{j=0}^{d_a-1}\sum_{i=1}^{\rho_{a,j}}\nu_i^{(a,j)}}$, respectively; finally, consider the complete sample matrices collecting all nodes for different $a\in\{1,\dots,N\}$
\begin{align*}
  	\Sigma_k^{n_1d_1+\cdots +n_Nd_N}\coloneq&\big((\Sigma_k^{n_1d_1+\cdots +n_Nd_N})_1,  ,\dots,(\Sigma_k^{n_1d_1+\cdots +n_Nd_N})_N \big),\\
	\Sigma_{[k,n_1d_1+\cdots +n_Nd_N]}\coloneq&\big((\Sigma_{[k,n_1d_1+\cdots +n_Nd_N]})_1,\dots,(\Sigma_{[k,n_1d_1+\cdots +n_Nd_N]})_N\big).
\end{align*}
\end{definition}
We use the word partial in the sense that they are linked to one of the involved prime polynomials.
We now proceed as we have done in previous situations just changing nodes and sample matrices as we have indicated. Then,
\begin{definition}
	We say that $\mathcal {N}_{k,n_1d_1+\cdots +n_Nd_N}$  is a poised set if the sample matrix is non singular
	\begin{align*}
	\det\Sigma_k^{n_1d_1+\cdots +n_Nd_N}\neq 0.
	\end{align*}
\end{definition}

\begin{theorem}\label{theorem:the big deal}
	For a poised set of nodes $\mathcal {N}_{k,n_1d_1+\cdots +n_Nd_N}$ in the algebraic hypersurface $\bigcup\limits_{a=1}^N Z(\mathcal Q_a)$  the transformed orthogonal polynomials can be expressed in terms of the original ones as the following last quasi-determinantal expression
	\begin{align*}
	TP_{[k]}(\x)=  \frac{\big(\prod_{a=1}^N(\mathcal R_a(\boldsymbol\Lambda))^{d_a}\big)_{[k],[k+n_1d_1+\cdots +n_Nd_N]}}{\mathcal R_1(\x)^{d_1}\cdots \mathcal R_N(\x)^{d_N}}
	\Theta_*\PARENS{\begin{array}{c|c}
		\Sigma^{n_1d_1+\cdots +n_Nd_N}_k &\begin{matrix}
		P_{[k]}(\x) \\ \vdots \\ P_{[k+n_1d_1+\cdots +n_Nd_N-1]}(\x)
		\end{matrix}\\\hline
		\Sigma_{[k,n_1d_1+\cdots +n_Nd_N]} & P_{[k+n_1d_1+\cdots +n_Nd_N]}(\x)
		\end{array}}.
	\end{align*}
\end{theorem}
\begin{proof}
	Proposition \ref{resolventd} gives
	\begin{align*}
	\omega_{[k],[k+n_1d_1+\cdots +n_Nd_N]}\Sigma_{[k,n_1d_1+\cdots +n_Nd_N]}+
	\big(\omega_{[k],[k]},\cdots,	\omega_{[k],[k+n_1d_1+\cdots +n_Nd_N-1]}\big)\Sigma_k^{n_1d_1+\cdots +n_Nd_N}=0,
	\end{align*}
so that
	\begin{multline*}
	(\omega_{[k],[k]},\dots,\omega_{[k],[k+n_1d_1+\cdots +n_Nd_N-1]})\\=-\big(\prod_{a=1}^N(\mathcal R_a(\boldsymbol\Lambda))^{d_a})\big)_{[k],[k+n_1d_1+\cdots +n_Nd_N]}\Sigma_{[k,n_1d_1+\cdots +n_Nd_N]}\big(\Sigma_k^{n_1d_1+\cdots +n_Nd_N}\big)^{-1}
	\end{multline*}
	and  $\mathcal R_1(\x)^{d_1}\cdots \mathcal R_N(\x)^{d_N}TP(\x)=\omega P(\x)$ gives the result.
\end{proof}
As in previous discussion we shift from the field $\R$ to its algebraic closure $\C$, laying the algebraic hypersurfaces and nodes in the $D$-dimensional complex affine space.
\begin{definition}
	We  introduce the partial derived Vandermonde matrices 
	\begin{align*}
	(\mathcal V_{k}^{n_1d_1+\cdots +n_Nd_N})^{(a,j)}_i\coloneq&
	\Big(\dfrac{\partial^j \chi^{[k+n_1d_1+\cdots +n_Nd_N]}}{\partial \n^{(a,j)}_i}(\boldsymbol p^{(a,j)}_{i,1}) ,\dots,
	\dfrac{\partial^j \chi^{[k+n_1d_1+\cdots +n_Nd_N]}}{\partial \n^{(a,j)}_i}(\boldsymbol p^{(a,j)}_{i,\nu^{(a,j)}_i})\Big),
	\end{align*}
	that belong to  $\C^{N_{k+n_1d_1+\cdots +n_Nd_N-1}\times \nu^{(a,j)}_i}$. For $ j\in\{0,\dots,d_a-1\}$ and $i\in\{1,\dots,\rho_{a,j}\}$,
	the deriveºd  Vandermonde matrix is
	\begin{align*}
	(\mathcal V_{k}^{n_1d_1+\cdots +n_Nd_N})^{(a,j)}\coloneq&
	\big((\mathcal V_{k}^{n_1d_1+\cdots +n_Nd_N})^{(a,j)}_1,\dots,(\mathcal V_{k}^{n_1d_1+\cdots +n_Nd_N})^{(a,j)}_{\rho_{a,j}}\big)\in\C^{N_{k+n_1d_1+\cdots +n_Nd_N-1}\times \sum_{i=1}^{\rho_{a,j}}\nu^{(a,j)}_i},
	\end{align*}
	and the partial multivariant confluent Vandermonde matrix in $\C^{N_{k+n_1d_1+\cdots +n_Nd_N-1}\times \sum_{j=0}^{d_a-1}\sum_{i=1}^{\rho_{a,j}}\nu^{(a,j)}_i},$
	\begin{align*}
	(\mathcal V^{n_1d_1+\cdots +n_Nd_N}_k)_a\coloneq \big( (\mathcal V_{k}^{n_1d_1+\cdots +n_Nd_N})^{(a,0)},(\mathcal V_{k}^{n_1d_1+\cdots +n_Nd_N})^{(a,1)},\dots,(\mathcal V_{k}^{n_1d_1+\cdots +n_Nd_N})^{(a,d_a-1)}\big).
	\end{align*}
	Finally, the complete confluent Vandermonde matrix in $\in\C^{N_{k+n_1d_1+\cdots +n_Nd_N-1}\times r_{k,n_1d_1+\cdots +n_Nd_N}},$
		\begin{align*}
	\mathcal V^{n_1d_1+\cdots +n_Nd_N}_k\coloneq \Big((\mathcal V^{n_1d_1+\cdots +n_Nd_N}_k)_1,\dots,(\mathcal V^{n_1d_1+\cdots +n_Nd_N}_k)_N\Big).
		\end{align*}
\end{definition}
Again we find the factorization
$\Sigma_k^{n_1d_1+\cdots +n_Nd_N}=S_k^{n_1d_1+\cdots +n_Nd_N}
\mathcal V_{k}^{n_1d_1+\cdots +n_Nd_N}$,
with $S^m_k$ as in \eqref{eq:S-slice},
so that  $\operatorname{Ker}\mathcal V^{n_1d_1+\cdots +n_Nd_N}_k \subset	\operatorname{Ker} \Sigma_k^{n_1d_1+\cdots +n_Nd_N}$.
For $\mathcal N_k^{n_1d_1+\cdots +n_Nd_N}$ to be  poised  we must request to  $\mathcal V_k^{n_1d_1+\cdots +n_Nd_N}$ to have full column rank matrix:
$\dim\operatorname{Im}\mathcal V_k^{n_1d_1+\cdots +n_Nd_N}=r_{k,n_1d_1+\cdots +n_Nd_N}$. Again, Theorems \ref{teorema:ideal-vandermonde} and \ref{theorem:poised-fullrank} and Propositions  \ref{pro:matrix-ideal} and \ref{pro:idel-vandermonde2} hold true for our polynomial $\mathcal Q=\mathcal R^d$ and the corresponding confluent Vandermonde matrix $\mathcal V^m_k$, $m={n_1d_1+\dots+n_Nd_N}$.

\begin{theorem}\label{HH}
%
		The node set $\mathcal N_{k,n_1d_1+\cdots +n_Nd_N}\subset\C^D$ is    poised if it does not exist a polynomial  $V\neq \prod_{a=1}^N\mathcal R_a^{d_a}$, $\deg V\leq k+n_1d_1+\cdots +n_Nd_N-1$, such that 
		$\mathcal N_i^{(a,j)}\subset Z\Big(\dfrac{\partial^j V}{\partial \n^{(a,j)}_i}\Big)$, 	for $a=1,\dots,N$,  $j\in\{1,\dots,d_a-1\}$ and $i\in\{1,\dots,\rho_{a,j}\}$.
\end{theorem}
\begin{proof}
A vector $v=(v_i)_{i=1}^{N_{k+n_1d_1+\cdots +n_Nd_N-1}}\in(\operatorname{Im}\mathcal V_k^{n_1d_1+\cdots+n_Nd_N})^\perp$ if for the corresponding polynomial $V=\sum_{l=1}^{N_{k+n_1d_1+\cdots +n_Nd_N-1}}v_l\x^{\q_l}$  the polynomials $\dfrac{\partial^j V}{\partial\n^{(a,j)}_i}$ cancel at $\mathcal N^{(a,j)}_i$. Notice  that the polynomial
	$V= \x^\q\prod_{b=1}^N\mathcal R_b^{d_b}$ is  such $\dfrac{\partial^j V}{\partial\n^{(a,j)}_i}$ do cancel  at $\cup_{b=1}^N Z(\mathcal R_b)$ for $j=0$ and also at $Z(\mathcal R_a)$ for $j\in\{1,\dots,d_a-1\}$ and $i\in\{1,\dots,\rho_{a,j}\}$. Hence, 
	\begin{align*}
	\Big(\prod_{a=1}^N\mathcal R_a^{d_a}\Big)_{k+n_1d_1+\cdots +n_Nd_N-1}\subseteq (\operatorname{Im}\mathcal V_k^{n_1d_1+\cdots +n_Nd_N})^\perp,
	\end{align*}
	which considered at the light of the condition  $\dim (\prod_{a=1}^N\mathcal R_a^{d_a})_{k+n_1d_1+\cdots +n_Nd_N-1}=N_{k-1}$ implies that no further constraint can be allowed or the full column rank would not be achievable. 
	Thus, the set  is poised if  $\operatorname{Im}\mathcal V^{n_1d_1+\cdots +n_Nd_N}_k=(\mathcal R^d)_{k+n_1d_1+\cdots +n_Nd_N-1}^\perp$.
\end{proof}
\begin{pro}
In order to have a  poised set the  polynomial $\dfrac{\partial^j (\mathcal R_1^{p_1}\cdots\mathcal R_N^{p_N})}{\partial\n_i^{(a,j)}}$  can not cancel at $Z(\mathcal R_1\cdots\mathcal R_N)$ for $0\leq p_1<d_1,\cdots,  0\leq p_N<d_N$.
Moreover, 	when the set of nodes is  poised we can ensure the following bounds for the node subset cardinals   
	\begin{align}\label{eq:lower bound general}
	|\mathcal N_a|&\geq \Big\lceil\frac{k}{d_a}\Big\rceil+n_a,\\
	|\mathcal N^{(a,0)}|&\leq r_{k+n_1d_1+\cdots+n_Nd_N-n_a,n_a}\label{eq:partial upper bound 1},\\
	|\mathcal N^{(a,j)}_i|&\leq r_{k+n_1d_1+\cdots+n_Nd_N-d^{(a,j)}_i,n_a}\label{eq:partial upper bound 2}.
	\end{align}
Here 
\begin{align*}
d^{(a,j)}_i\coloneq\deg \frac{\partial^j(\mathcal R_a^{d_a})}{\partial \n^{(a,j)}_i}
\end{align*}
and the function $\lceil x\rceil$  gives the smallest integer $\geq x$.
\end{pro}
\begin{proof}
	If $\dfrac{\partial^j (\mathcal R_1^{p_1}\cdots\mathcal R_N^{p_N})}{\partial\n_i^{(j)}}$  cancels at $Z(\mathcal R)$ then 
	\begin{align*}
	\dim \big((\operatorname{Im}\mathcal V_k^{n_1d_1+\cdots+n_Nd_N})^\perp\big)\geq N_{k+n_1d_1+\cdots+n_Nd_N-(p_1n_1+\cdots+p_Nd_N)}
	\end{align*}
	 and
	the set is not  poised. 
	
To prove \eqref{eq:lower bound general}  let us consider first order polynomials
 $\pi^{(a,j)}_{i,l}(\z)$   which cancels at $\boldsymbol p^{(a,j)}_{i,l}$,  $\pi^{(a,j)}_{i,l}(\boldsymbol p^{(a,j)}_{i,l})=0$, and construct the polynomial $\Pi_a=\Big(\prod_{\substack{j=1,\dots,d_a\\i=1,\dots,\rho_{a,j}\\l=1,\dots,\nu_i^{(a,j)}}}\pi^{(a,j)}_{i,l}\Big)$, $\deg \Pi_a=|\mathcal N_a|$. Then,  
the polynomial $V=\mathcal R_1^{d_1}\cdots\mathcal R_{a-1}^{d_{a-1}}\mathcal R_{a+1}^{d_{a+1}}\cdots\mathcal R_N^{d_N}\Pi_a^{d_a}$, $\deg V=n_1d_1+\dots+n_Nd_N+(|\mathcal N_a|-n_a)d_a$, has  the nodes  among its zeroes $\mathcal N^{n_1d_1+\cdots+n_Nd_N}_k\subset Z(V)$ and  $\dfrac{\partial^j V}{\partial\n^{(b,j)}_i}$ do cancel at $\mathcal N^{(b,j)}_i$ for all $b\in\{1,\dots,N\}$. Thus, we should request 
$n_1d_1+\dots+n_Nd_N+(|\mathcal N_a|-n_a)d_a\geq k+n_1d_1+\cdots +n_Nd_N$; i.e., $(|\mathcal N_a|-n_a)d_a\geq k$ and the result follows.

For \eqref{eq:partial upper bound 1} observe that all the columns  in the  Vandermonde block $(\mathcal V_k^{n_1d_1+\cdots+n_Nd_N})^{(a,0)}$  imply no directional partial derivatives and are evaluated at nodes which belong to $Z(\mathcal R_a)$. Therefore, 
\begin{align*}
(\mathcal R_a)_{k+n_1d_1+\cdots+n_Nd_N-1}^\perp\supseteq\operatorname{Im}(\mathcal V_k^{n_1d_1+\cdots+n_Nd_N})^{(0)}_a
\end{align*}
and the maximum rank achievable for this block is $N_{k+n_1d_1+\cdots+n_Nd_N-1}-N_{k+n_1d_1+\cdots+n_Nd_N-n_a-1}=r_{k+n_1d_1+\cdots+n_Nd_N-n_a,n_a}$. The columns of the block 
$(\mathcal V_k^{n_1d_1+\cdots+n_Nd_N})^{(a,j)}_i$ are vinculated to 
$Z\big(\frac{\partial^j(\mathcal R_a^{d_a})}{\partial \n^{(a,j)}_i}\big)$
and consequently  $\big(\frac{\partial^j(\mathcal R_a^{d_a})}{\partial \n^{(j)}_i}\big)_{k+k+n_1d_1+\cdots+n_Nd_N-1}^\perp\supseteq \Big(\operatorname{Im}(\mathcal V_k^{k+n_1d_1+\cdots+n_Nd_N})^{(j)}_i\Big)$; hence, the maximum possible rank for this block  is $N_{k+k+n_1d_1+\cdots+n_Nd_N-1}-N_{k+k+n_1d_1+\cdots+n_Nd_N-1-d^{(a,j)}_i}$.
\end{proof}

\end{document}